\documentclass[12pt,leqno,twoside]{article} 
\input{twisted.sty}
\usepackage{bezier,eufrak}
\begin{document}
\title{On representations of twisted group rings}
\author{Matthias K\"unzer}
\maketitle

\begin{small}
\begin{quote}
\begin{center}{\bf Abstract}\end{center}\vspace*{2mm}
We generalize certain parts of the theory of group rings to the twisted case. Let $G$ be a finite group acting (possibly trivially) on a field $L$ of characteristic
coprime to the order of the kernel of this operation. Let $K\tm L$ be the fixed field of this operation, let $S$ be a discrete valuation ring with field of fractions $K$, 
maximal ideal generated by $\pi$ and integral closure $T$ in $L$. We compute the colength
of $T\wr G$ in a maximal order in $L\wr G$. Moreover, if $S/\pi S$ is finite, we compute the $S/\pi S$-dimension of the center of $T\wr G/\Jac(T\wr G)$. If this quotient is
split semisimple, this yields a formula for the number of simple $T\wr G$-modules, generalizing Brauer's formula.
\end{quote}
\end{small}

\renewcommand{\thefootnote}{\fnsymbol{footnote}}
\renewcommand{\thefootnote}{\arabic{footnote}}

\begin{footnotesize}
\renewcommand{\baselinestretch}{0.7}
\parskip0.0ex
\tableofcontents
\parskip1.2ex
\renewcommand{\baselinestretch}{1.0}
\end{footnotesize}

\setcounter{section}{-1}

\section{Introduction}

\subsection{Outline}
Let $G$ be a finite group acting on a field $L$ with fixed field $K$. The twisted 
group ring $L\wr G$ carries the multiplication $(\sigma y)(\tau z) = \sigma\tau y^\tau z$, where $\sigma,\tau\in G$, $y,z\in L$ 
(\footnote{This ring is also known as {\it skew group ring} of $G$ with coefficients in $L$, see e.g.\ {\bf\cite[\rm p.\ 59]{Ka90}}. Using the terminology {\it twisted group ring}
of {\bf\cite[\rm\S 28]{CR81}}, we should mention that there are more general twisted group rings still, involving $2$-cocycles of $G$ with coefficients in $L^\ast$.}). 

We denote by $N$ the kernel of the operation of $G$ on $L$.
If $N = G$, then $K = L$ and $L\wr G$ is nothing but the untwisted classical group ring $KG$. If $N = 1$, then $L\wr G$ has
the only simple module $L$, and is isomorphic to $\End_K L$. 

If the characteristic of $L$ and the order of $N$ are coprime, then $L\wr G$ is semisimple, and even separable over $K$. Many of the
apparent difficulties we shall encounter vanish if the characteristic of $L$ and the order of $G$ are coprime. For instance, in this case the Plancherel
formula (\ref{PropI_3_5}) is immediate.

The Wedderburn isomorphism sends an element of $L\wr G$ to the tuple of its operations on the simple modules, such an operation being considered as an endomorphism over the 
respective $L\wr G$-linear endomorphism skew field. The Plancherel formula yields a Fourier inversion formula for this isomorphism. Schur relations may then be
deduced from composition of the Wedderburn isomorphism with its Fourier inversion. A theory of characters and a Brauer-Nesbitt theorem ensue.

Let $S\tm K$ be a discrete valuation ring with 
field of fractions $K$, maximal ideal generated by $\pi$ and residue field $\b S := S/\pi S$ of characteristic $p \geq 0$. The $\pi$-adic valuation of an element $y\in K$ is 
denoted by $v(y)$. Let $T\tm L$ be the integral closure of $S$ in $L$, and denote $\b T := T/\pi T$. 
Choosing $T\wr G$-lattices inside the simple $L\wr G$-modules, we may restrict the Wedderburn isomorphism to a full embedding of $T\wr G$ into a tuple of matrix rings over
certain extensions of $S$. We calculate the colength of this embedding, that is, the $S$-linear length of its cokernel.

Now suppose $\b S$ to be finite.
The question for the number of indecomposable projective \mb{$T\wr G$-modules} leads to the question for the number of simple $\b T\wr G$-modules. If \mb{$\b T\wr G/\Jac(\b T\wr G)$} 
is a split semisimple $\b S$-algebra, then 
this number coincides with the $\b S$-linear dimension of $Z(\b T\wr G/\Jac(\b T\wr G))$. This dimension in turn is given by
a sum of certain vectorspace dimensions, indexed by the $p'$-classes of $G$. In the untwisted case, each of these dimensions equals $1$. Thus we recover a theorem of {\sc Brauer.}

\subsection{Motivation}

The untwisted group ring $SG$ has the principal module $S$, endowed with the trivial $G$-operation. This module is projective if and only if the order of $G$ is invertible in $S$. 

The twisted group ring $T\wr G$ has the principal module $T$, endowed with the Galois operation of $G$. This module is projective if and only if the order of $N$ is 
invertible in $S$ and the extension $T/S$ is tamely ramified ({\sc Speiser,} {\sc E.\ Noether,} cf.\ \ref{RemI3_2}).

Therefore, wildness of $T/S$ corresponds to a nontrivial homological behaviour of $T$ over $T\wr G$, which one is led to study. This seems to be more easily accessible
than the related homological behaviour of $T$ over $SG$. For instance, we have $\Ext_{T\wr G}^\ast(T,T) \iso \Ext_{SG}^\ast(S,T)$, and in particular cases it is possible 
to calculate the left hand side in order to obtain the right hand side. Here, we will restrict our considerations to some general properties of twisted group rings, however.

\subsection{Known Results}

Beautiful results on the classical case $N = G$ abound. It might be worthwhile to attempt to extend certain of them to the twisted case. For example, so far there is no theory on 
induction and restriction that allows to independently vary the group and the ring (cf.\ {\bf\cite[\rm ch.\ 4]{Ka90}}, {\bf\cite[\rm\S\S\ 10, 11, 19, 20]{CR81}})
(\footnote{For this purpose, it might be useful to consider a larger class of rings with group operation than merely Galois extensions of fields -- cf.\ section \ref{SubsecMor}.}).

On the other hand, there are several results for the case $N = 1$. The theorem of {\sc Auslander-Goldman} states that $T\wr G$ is isomorphic to $(S)_{|G|}$ if and only if
$T/S$ is unramified and all occurring residue field extensions are separable {\bf\cite[\rm 28.5]{CR81}}. The theorem of {\sc Auslander-Rim} 
states that $T\wr G$ is hereditary if and only if $T/S$ is tamely ramified {\bf\cite[\rm 28.7]{CR81}}. 

If we admit wild ramification, or if nontrivial $2$-cocycles enter, the situation gets more involved. Still, we assume $N = 1$.
\begin{itemize}
\item {\sc Wilson}'s theorem calculates the locally free class group of $\Ol_L\wr G$ to be the class group of a maximal $\Ol_K$-order containing it, provided $\Ol_L/\Ol_K$
is a finite extension of Dedekind domains giving rise to an extension $L/K$ of global fields with Galois group $G$ {\bf\cite[\rm 50.64, 49.33, 49.32]{CR87}}. 
\item {\sc Williamson} gives a criterion for $\pi$-principal heredity of $T\wr_f G$, where $f$ is a $2$-cocycle of $G$ with values in $T^\ast$ {\bf\cite[\rm 3.5]{W66}}. 
In contrast to the above-mentioned theorems, it may occur that $T\wr_f G$ is hereditary, and even maximal, whereas $T/S$ is wildly ramified {\bf\cite[\rm 2.5]{W65}}. 
\item {\sc Benz} and {\sc Zassenhaus} showed that $T\wr_f G$ is contained in a unique minimal hereditary overorder, $S$ assumed to be complete {\bf\cite{BZ85}}.
\item {\sc Cliff} and {\sc Weiss} showed how invariants of this hereditary overorder depend on the inertia and ramification indices of $T/S$, and
the Schur index of $L\wr_f G$ {\bf\cite{CW85}}.
\item If $[L:K] = 2$, {\sc Chalatsis} and {\sc Theohari-Apostolidi} classified the indecomposable $\Ol_L\wr G$-lattices {\bf\cite[\rm Th.\ 2]{CTA86}}.
\item {\sc Weber} described the image of the Wedderburn embedding for $T = \Z_{(p)}[\zeta_{p^m}]$, $G = C_p$ and $S = \Z_{(p)}[\zeta_{p^{m-1}}]$ {\bf\cite{We00}}.
\end{itemize}

For $N$ arbitrary, {\sc Nakayama} and {\sc Shoda} investigated twisted group rings $L\wr G$ over a field $L$ such that the order of $G$ is invertible in $L$, relating
the representations of $L\wr G$ and $KN$ {\bf\cite{NaSh35}}. Omitting this invertibility condition, {\sc Bonami} investigated centers and radicals in {\bf\cite{Bo84}}.

Moreover, the ring $K\wr_f G$ has been studied thoroughly, where $f$ is a $2$-cocycle with coefficients in the trivial $G$-module $K^\ast$; one reason being its connection
to Clifford theory (cf.\ {\bf\cite[\rm 11.19, 11.20]{CR81}}). See also {\bf\cite[\rm (10.5)]{T95}}.

\subsection{Results}

Let $g := |G|$, $n := |N|$ and $h := [L:K]$, so $nh = g$. Let $\{X_i\;|\; i\in [1,k]\}$ be a complete set of nonisomorphic simple $L\wr G$-right modules and let 
$K_i := \End_{L\wr G} X_i$. Let $x_i$ denote the dimension of $X_i$ as a right module over $K_i$, let $r_i$ denote the degree $[K_i:K]$ and let $d_i^2$ denote the degree of $K_i$ 
over its center. We have $x_i d_i/h\in\Z$ (\ref{LemI_3_2}). Let $V_i\tm X_i$ be a finitely generated $T\wr G$-module such that $KV_i = X_i$ and such that $V_i$ is free as a module over 
$S_i := \End_{T\wr G} V_i$. 

The colength of an embedding of $S$-modules is defined to be the Jordan-H\"older multiplicity of $\b S$ in its quotient. If $\b S$ is finite, then 
$\mb{colength} = \log_{|\b S|}(\mb{index})$. 

Let $S_i^+ := \{y\in K_i\; |\; \Tr_{K_i/K}(yS_i)\tm S\}$ and let $\delr_{S_i/S}$ be the colength of $S_i\tm S_i^+$. 

Let $T^+ := \{y\in L\; |\; \Tr_{L/K}(yT)\tm S\}$ and let $\delr_{T/S}$ be the colength of $T\tm T^+$. 

{\bf Theorem} (\ref{PropI7}). {\it The colength of the Wedderburn embedding
\[
T\wr G \;\;\hraa{\omega}\; \prodd{i\in [1,k]} \End_{S_i} V_i
\]
is given by
\[
\frac{1}{2}\left(g(\delr_{T/S} + v(n)h) - \sumd{i\in [1,k]} x_i^2 (\delta_{S_i/S} + v(x_i d_i/h) r_i ) \right).
\]
In particular, this integer is positive or zero.
}

The essential ingredient for (\ref{PropI7}) is the Plancherel formula (\ref{PropI_3_5}). Using an analogous central Plancherel formula (\ref{PropCB3}), we obtain a 
central colength formula (\ref{PropZ4}).

Schur relations allow to derive a Brauer-Nesbitt theorem. Let $t$ be the maximal elementary divisor exponent of the Gram matrix of the trace bilinear 
form on $T$ induced by $\Tr_{L/K}$. Note that $v(x_i) \leq v(g) + t - v(d_i)$ (\ref{LemI15_1}).

{\bf Theorem} (\ref{CorBN1}). {\it Let $i\in [1,k]$. Suppose that $K = K_i$ and that $v(x_i) = v(g) + t$. Then the reduction $V_i/\pi V_i$ of the simple $\Lambda$-lattice $V_i$
is a simple $\b T\wr G$-module.}

Moreover, we shall count simple modules.

{\bf Theorem} (\ref{CorI17}, \ref{ThBrauer}).
\vspace*{-3mm}
\begin{itemize}
\item[(i)] {\it The $K$-linear dimension of the center $Z(L\,\wr\, G)$ equals the number of conjugacy classes of $N$. In particular,
if $L\wr G$ is split semisimple over $K$, then this number equals the number of isoclasses of simple $L\wr G$-modules.}

\item[(ii)]
{\it Let $\b T_0 := T/\Jac(T)$. Let
\[
V_\sigma \; :=\; \b S\spi{(y^\sigma - y)z,\; y^\rho - y\; |\; y,z\in\b T_0,\; \rho\in C_G(\sigma)}\; \tm\; \b T_0\; .
\]
Let $\ClGpp$ be a set of representatives of $p'$-classes of $G$. We have
\[
\dim_{\b S} Z((T\wr G)/\Jac(T\wr G))\; =\; \sum_{\sigma\in\sClGpp}\dim_{\b S}(\b T_0/V_\sigma)\; .
\]
In particular, if $T\wr G/\Jac(T\wr G)$ is split semisimple over $\b S$, then this number equals the number of isoclasses of simple $T\wr G$-modules, or, which is the same, 
the number of isoclasses of indecomposable projective $T\wr G$-modules.}
\end{itemize}

If $\b T_0$ is a field, and if $\sigma$ is not contained in the kernel $N_0$ of the operation of $G$ on $\b T_0$, 
then $\b T_0/V_\sigma = 0$. This indicates a similarity of (ii) to (i). 

\subsection{Acknowledgements}

I would like to thank {\sc G.\ Nebe} for pointing out {\sc Plesken'}s method to calculate colengths, and for several helpful discussions. I would like to thank the referee for
providing the isomorphism (\ref{LemMor1}), yielding a splitting field for $L\wr G$ (\ref{CorMor3}).

Dedicated to {\sc K.\ W.\ Roggenkamp} on the occasion of his 60th birthday.

\subsection{Conventions}

\begin{footnotesize}
\begin{itemize}
\item[(i)] Composition of maps is written on the right, $\lraa{a}\lraa{b} = \lraa{ab}$. Exception is made for `standard' maps, such as traces, characters, Frobenius \dots
\item[(ii)] If $A$ is an assertion, which might be true or false, we let $\{ A\} = 1$ if $A$ is true, and $\{ A\} = 0$ if $A$ is false.
\item[(iii)] For $a,b\in\Z$, we denote by $[a,b] := \{c\in\Z\;|\; a\leq c\leq b\}$ the integral interval.
\item[(iv)] Given $a\in\Z$ and a prime $p$, we denote by $a[p] := p^{v_p(a)}$ the $p$-part of $a$.
\item[(v)] Given elements $x, y$ of some set $X$, we let $\dell_{x,y} = 1$ in case $x = y$ and $\dell_{x,y} = 0$ in case $x\neq y$.
\item[(vi)] Unless mentioned otherwise, modules are finitely generated right modules.
\item[(vii)] Let $A$ be a commutative ring. The trace of an element of an $A$-free $A$-algebra is the trace of the $A$-linear map given by right multiplication with this element. 
\item[(viii)] Given a commutative ring $A$ and $x,y,z\in A$, we write $x\con_z y$ if $x - y \in Az$.
\item[(ix)] Given a ring $A$ and $n\geq 1$, we denote by $(A)_n$ the ring of $n\ti n$-matrices over $A$.
\item[(x)] Given a ring $A$, we denote by $\Jac(A) \tm A$ its Jacobson radical. 
\item[(xi)] Given an element $\sigma$ of a finite group $G$, we let $o(\sigma)$ denote its order.
\end{itemize}
\end{footnotesize}        
\section{Rational}
\label{SubsecAss}

\subsection{The Plancherel formula}

\begin{Definition}
\label{DefI0}\rm
Let $C$ be a commutative ring, acted upon by a group $G$ with fixed ring $B := \mb{Fix}_G C$, an element $\tau\in G$ sending an element $y\in C$ to $y^\tau$.  
The {\it twisted group ring} $C\wr G$ is defined as follows. As a right $C$-module, it is free on the underlying set of the group $G$. The multiplication is defined by
\[
\left(\sumd{\sigma \in G}\sigma y_\sigma\right)\left(\sumd{\tau \in G}\tau z_\tau\right) 
\; :=\;  \sumd{\rho\in G}\rho \left(\sumd{\tau\in G} y_{\rho\tau^{-1}}^\tau z_\tau\right),
\]
where $y_\sigma, z_\tau\in C$ for $\sigma,\tau\in G$. Or, less formally, we extend $(\sigma y)(\tau z) = (\sigma\tau) (y^\tau z)$ $B$-linearly. Note that $C\wr G$ is a $B$-algebra.
\end{Definition}

\begin{Notation}
\label{SetupI0_5}\rm
Let $G$ be a finite group acting on a field $L$ with fixed field $K := \mb{Fix}_G L$. Let $H$ be the image of the operation map $G\lra\Aut_{\!\mb{\scr fields}}\, L$, i.e.\ 
$H = \Gal(L/K)$, and let $N$ denote the kernel of this operation map. The orders are denoted by
\[
\begin{array}{rcl}
h & := & |H| \; =\; [L:K] \\
g & := & |G|  \\ 
n & := & |N|\; , \\ 
\end{array}
\]
so $g = hn$. Let $u\in L$ be such that $\Tr_{L/K}(u) = 1$. Choose a $K$-basis $(y_l)_{l\in [1,h]}$ of $L$, and
let $(y^\ast_l)_{l\in [1,h]}$ be its dual basis with respect to $\Tr_{L/K}$, i.e.\ $\Tr_{L/K}(y_l y^\ast_m) = \dell_{l,m}$ for $l,m\in [1,h]$. Whenever useful, we assume $y_1 = 1$.
\end{Notation}

\begin{Lemma}[{Maschke, cf.\ {\bf\cite[\rm 28.7]{CR81}}}]
\label{LemI2}
The $K$-algebra $L\wr G$ is semisimple if and only if $n$ is invertible in $K$. 

\rm
Suppose $n$ to be invertible in $K$. Given an epimorphism of $L\wr G$-modules $X\lraa{f} X''$, we choose an $L$-linear coretraction $X''\lraa{i_0} X$ and let
\[
\begin{array}{rcl}
X'' & \lraa{i} & X \\
x'' & \lramaps & n^{-1}\sumd{\sigma\in G} (x''\sigma)i_0\cdot u\sigma^{-1}\; , \\
\end{array}
\]
which is $L\wr G$-linear and satisfies $if = 1_{X''}$.

Conversely, suppose $L\wr G$ to be semisimple. Since $L\wr G\lra L : \sigma y\lramaps y$ is an epimorphism, it has a coretraction, which is necessarily of the form 
$L\lra L\wr G : 1\lramaps z\sum_{\sigma\in G} \sigma$ for some $z\in L$. Composition yields $n\cdot\Tr_{L/K}(z) = 1$.
\end{Lemma}

Henceforth, we \fbox{assume $n$ to be invertible in $K$\rule[-1.3mm]{0mm}{5mm}.}
Now, alternatively, semisimplicity ensues from the following lemma {\bf\cite[\rm 7.18, 7.20]{R75}}

\begin{Lemma}
\label{LemI2_5}
The $K$-algebra $L\wr G$ is separable.

\rm
We have to show that the $(L\wr G)^o\ts_K (L\wr G)$-linear epimorphism 
\[
\begin{array}{rllcl}
(L\wr G)^\mb{\scr o} & \ts_K & (L\wr G) & \lra     & L\wr G \\ 
\sigma u             & \ts   & v\rho    & \lramaps & \sigma u v\rho \\
\end{array}
\]
is split, where $\sigma,\rho\in G$, $u,v\in L$. Let $\sum_{i\in I} \xi_i\ts\eta_i\in L\ts_K L\iso\prod_H L$ be the element that 
satisfies $\sum_{i\in I}\xi_i\eta_i^\lambda = \dell_{1,\lambda} u$ for $\lambda\in H$. Then
\[
\begin{array}{rcrll}
L\wr G & \lra     & (L\wr G)^\mb{\scr o}                                 & \ts_K & (L\wr G) \\ 
1      & \lramaps & n^{-1}\sumd{\kappa\in G,\; i\in I} \kappa\xi_i\!\!\! & \ts   & \!\!\!\!\!\!\!\eta_i\kappa^{-1} \\
\end{array}
\]
is a coretraction as sought for.
\end{Lemma}

\begin{Lemma}[Dedekind]
\label{LemI1}
Equipped with the operation $y\cdot (\sigma z) := y^\sigma z$, where $y, z\in L$, $\sigma\in G$, the field $L$ becomes an absolutely simple module over $L\wr G$, called the 
{\rm principal module.} We have $\End_{L\wr G} L = K$. The module $L$ is the simple module of the block $L\wr H$ of $L\wr G$, unique up to isomorphism.

\rm
Consider the map
\[
\begin{array}{rcl}
L\wr H   & \lraa{\omega_H} & \End_K L \\
\sigma y & \lramaps        & (z\lramaps z^\sigma y). \\
\end{array}
\]
Both source and target are of dimension $h^2$ over $K$. Injectivity of $\omega_H$ is the content of Dedekind's Lemma.
\end{Lemma}

\begin{Definition}
\label{DefI2_6}\rm
Let
\[
\fbox{$\;\;
(\sigma y,\tau z) \; :=\; \dell_{\sigma,\tau^{-1}}\Tr_{L/K}(y^\tau z) = \frac{1}{n} \Tr_K(u(-)\sigma y\tau z) 
\;\;$}
\]
define a $K$-bilinear form on $L\wr G$, where $\sigma,\,\tau\in G$, $y,\, z\in L$, and where $\xi (-)\xi'$ denotes the $K$-linear endomorphism of $L\wr G$ that sends $\eta$ to 
$\xi\eta\xi'$, where $\xi,\xi'\in L\wr G$. This bilinear form is symmetric, associative and nondegenerate. 
\end{Definition}

\begin{Notation}
\label{SetupI2_7}\rm
Let $\{X_i\;|\; i\in [1,k]\}$ be a complete set of nonisomorphic simple $L\wr G$-right modules and let $K_i := \End_{L\wr G} X_i$ be the respective endomorphism skew field. Let
$x_i$ denote the dimension of $X_i$ as a {\it left} module over $K_i$, let $r_i$ denote the degree $[K_i:K]$. Choose $X_1 = L$, whence $K_1 = K$. Choose a maximal commutative
subfield $E_i$ of $K_i$ that is separable over $K$, and let $c_i := [Z(K_i):K]$ and $d_i := [E_i : Z(K_i)] = [K_i : E_i]$ {\bf\cite[\rm 7.15]{R75}}. In particular, $r_i = c_i d_i^2$.
So altogether,
\[
(K\;\auf{c_i}{\tm}\; Z(K_i)\; \auf{d_i}{\tm}\; E_i\; \auf{d_i}{\tm}\; K_i) \; =\; (K\;\auf{r_i}{\tm} K_i)\; .
\]
Denote the operation of $L\wr G$ on $X_i$ by
\[
L\wr G\;\lraa{\omega_i}\;\End_{K_i} X_i \; .
\]
The Wedderburn isomorphism, of $K$-algebras, reads
\[
\begin{array}{rcl}
L\wr G   & \lraisoa{\omega} & \prodd{i\in [1,k]} \End_{K_i} X_i \\
\xi      & \lramaps         & (\xi\omega_i)_{i\in [1,k]}\; . \\
\end{array}
\]
In particular, $gh = \sum_{i\in [1,k]} r_i x_i^2$. Note that we may interpret $\omega_1$, up to isomorphism, as the surjective ring morphism $L\wr G\lra L\wr H$ induced by $G\lra H$.

Let $\tr_{K_i/Z(K_i)} : K_i\lra Z(K_i)$ denote the reduced trace {\bf\cite[\rm 9.6a, 9.3]{R75}}, and let 
\[
\tr_{K_i/K} \; :=\;  \Tr_{Z(K_i)/K}\circ \tr_{K_i/Z(K_i)} : K_i\lra K\; .
\]
\end{Notation}

\begin{Remark}
\label{RemI_3_0}\rm
The next aim is to provide the Plancherel formula (\ref{PropI_3_5}), which is needed in section \ref{SubsecResInt}. Note that if $h$ was assumed to be 
invertible in $L$, then the following considerations would simplify considerably, for then we would be able to write $(\sigma y,\tau z) = \frac{1}{g} \Tr_K((-)\sigma y\tau z)$. 
In our situation, however, we have to split up everything and to keep track of the various identifications. 
\end{Remark}

\begin{Notation}[auxiliary]
\label{NotI_3_1}\rm
Let $F$ be a finite galois extension of $K$ that contains $L$ and that contains $E_i$ for all $i\in [1,k]$. Given $i\in [1,k]$, we let $\Theta_i$ be the set of $K$-linear embeddings 
of $Z(K_i)$ into $F$, so in particular $|\Theta_i| = c_i$. For each $\theta\in \Theta_i$, we choose a prolongation $\h\theta: E_i\lra F$ of $\theta$.

By choice of a $K_i$-linear basis, we identify $X_i = K_i^{(x_i)}$, considered as a row of $\End_{K_i} X_i = (K_i)_{x_i}$. 

Given $\mu\in H$, we let $e_\mu\in F\ts_K L\iso\prod_{\lambda\in H} F$ be mapped to $(\dell_{\mu,\lambda})_\lambda$ under this isomorphism.
\end{Notation}

\begin{Remark}
\label{RemI_3_a}\rm
Suppose given $i\in [1,k]$. We have an isomorphism of $E_i$-algebras
\[
\begin{array}{rllrcl}
E_i & \ts_K & K_i & \lraiso  & \End_{E_i} K_i \\
e   & \ts   & w   & \lramaps & e(-)w  \\
\end{array}
\]
(cf.\ {\bf\cite[\rm pf.\ of 7.15, pf.\ of 7.11]{R75}}). Therefore, we have an isomorphism of $F$-algebras
\[
\begin{array}{rllrcl}
F & \ts_K & K_i & \lraiso  & \prod_{\theta\in\Theta_i} \End_F (F\liu{\h\theta\,}{\ts} K_i) \\
f & \ts   & w   & \lramaps & (f\liu{\h\theta\,}{(-)}w)_{\theta} \; , \\
\end{array}
\]
where we denote $F\liu{\h\theta\,}{\ts}K_i := F\ts_{E_i} K_i$ with respect to $E_i\hraa{\h\theta} F$, and where $f\liu{\h\theta\,}{(-)}w$ denotes the $F$-linear endomorphism
of $F\liu{\h\theta\,}{\ts} K_i$ that sends $f'\ts w'$ to $ff'\ts w'w$ (cf.\ {\bf\cite[\rm p.\ 100]{R75}}).
\end{Remark}

\begin{Lemma}
\label{LemI_3_4}
Given $i\in [1,k]$, $w\in K_i$ and $\theta\in\Theta_i$, we have 
\[
\sum_{\theta\in\Theta_i} \Tr_F(1\liu{\h\theta\,}{(-)}w) \; =\; \tr_{K_i/K}(w)\; \in\; K\; . 
\] 

\rm
Using an $E_i$-linear basis of $K_i$ as an $F$-linear basis of $F\liu{\h\theta\,}{\ts_K} K_i$, we see, denoting by $(-)w$ the $E_i$-linear multiplication endomorphism of $K_i$, that 
\[
\begin{array}{rcl}
\Tr_F(1\liu{\h\theta\,}{(-)}w) 
& = & (\Tr_{E_i}((-)w))\h\theta \\
& = & (\tr_{K_i/Z(K_i)} (w))\theta \\
\end{array}
\]
by definition of the reduced trace {\bf\cite[\rm 9.6a, 9.3, pf.\ of 7.15]{R75}}, cf.\ (\ref{RemI_3_a}).
\end{Lemma}

\begin{Lemma}
\label{LemI_3_b}
Representatives for the isoclasses of simple $F\ts_K L\wr G$-modules are given by $F\liu{\h\theta\,}{\ts} X_i$, where $i\in [1,k]$, $\theta\in\Theta_i$. Note that $X_i$ is a
left $E_i$-module by restriction from $K_i$. Given $f\in F$ and $\xi\in L\wr G$, $f\ts\xi$ acts on $F\liu{\h\theta\,}{\ts} X_i$ as
\[
(f\ts\xi)\omega_{i,\theta} \; :=\;  f\liu{\h\theta\,}{\ts} (\xi)\omega_i\; ,
\]
giving the component $F\ts_K L\wr G\;\lraa{\omega_{i,\theta}}\;\End_F F\liu{\h\theta\,}{\ts} X_i$ of the Wedderburn isomorphism.

\rm
A block of $(L\wr G)\omega$ splits under $F\ts_K -$ into
\[
\begin{array}{rcl}
F\ts_K \End_{K_i} X_i
& = & (F\ts_K K_i)_{x_i} \\
& \aufgl{(\ref{RemI_3_a})} & \prod_{\theta\in\Theta_i} (\End_F(F\liu{\h\theta\,}{\ts} K_i))_{x_i} \\
& = & \prod_{\theta\in\Theta_i} \End_F(F\liu{\h\theta\,}{\ts} X_i)\; . \\
\end{array}
\]
Now, $f\ts\phi$ corresponds to $(f\ts\phi_{a,b})_{a,b\in [1,x_i]}$, then to $((f\liu{\h\theta\,}{(-)}\phi_{a,b})_{a,b\in [1,x_i]})_\theta$, hence finally to 
$(f\liu{\h\theta\,}{\ts}\phi)_\theta$.
\end{Lemma}

\begin{Lemma}
\label{LemI_3_c}
Given $i\in [1,k]$, $\theta\in\Theta_i$ and $\phi\in\End_{K_i} X_i$, we have
\[
\sum_{\theta\in\Theta_i} \Tr_F(1\liu{\h\theta\,}{\ts} \phi)\; =\; \tr_{K_i/K} \Tr_{K_i}\phi \; .
\]

\rm
In fact,
\[
\begin{array}{rcl}
\sum_{\theta\in\Theta_i} \Tr_F(1\liu{\h\theta\,}{\ts} \phi)
& = & \sum_{a\in [1,x_i]}\sum_{\theta\in\Theta_i} \Tr_F (1\liu{\h\theta\,}{(-)}\phi_{a,a}) \\
& \aufgl{(\ref{LemI_3_4})} & \sum_{a\in [1,x_i]} \tr_{K_i/K} \phi_{a,a} \\
& = & \tr_{K_i/K} \Tr_{K_i} \phi\; . \\
\end{array}
\]
\end{Lemma}

\begin{Lemma}
\label{LemI_3_2}
Suppose given $i\in [1,k]$ and $\theta\in\Theta_i$. We have $x_i d_i/h\in\Z$. Moreover, given $z\in Z(K_i)\tm\End_{K_i} X_i$, we have
\[
\Tr_F\Big((e_\mu)\omega_{i,\theta}\Big) \; =\; x_i d_i/h\; .
\]

\rm
We decompose $F\liu{\h\theta\,}{\ts} X_i = \ds_{\lambda\in H} (F\liu{\h\theta\,}{\ts} X_i) e_\lambda$ as $F$-vectorspaces. Since 
$(F\liu{\h\theta\,}{\ts} X_i) e_\lambda\iso (F\liu{\h\theta\,}{\ts} X_i) e_{\sigma^{-1}\lambda}$ for 
$\sigma\in G$ via multiplication with $\sigma$, using $e_\mu \sigma = \sigma e_{\sigma^{-1}\mu}$, the $F$-linear dimension of $(F\liu{\h\theta\,}{\ts} X_i) e_\mu$ is given by 
$x_i d_i/h$. Now,
\[
\begin{array}{rcl}
\Tr_F((e_\mu)\omega_{i,\theta})
& = & \sum_{\lambda\in H} \Tr_F((e_\mu)\omega_{i,\theta}|_{(F\liu{\h\theta\,}{\ts} X_i) e_\lambda}) \\
& = & \Tr_F((e_\mu)\omega_{i,\theta}|_{(F\liu{\h\theta\,}{\ts} X_i) e_\mu})  \\
& = & x_i d_i/h\; . \\
\end{array}
\]
\end{Lemma}

\begin{Lemma}
\label{LemI_3_3}
Given $m\geq 1$ and $A,B\in (F)_m$, we have $\Tr_F(A(-)B) = \Tr_F(A)\Tr_F(B)$, where $A(-)B$ denotes the $F$-linear endomorphism of $(F)_m$ that sends $C\in (F)_m$ to $ACB$, and
where $\Tr_F(A)$ denotes the trace of $A$ as an element of $(F)_m$.
\end{Lemma}

\begin{Proposition}[Plancherel formula]
\label{PropI_3_5}
Given $\xi,\eta\in L\wr G$, we have
\[
(\xi,\eta)\; =\;\sum_{i\in [1,k]} \frac{x_i d_i}{g}\cdot \tr_{K_i/K} \Tr_{K_i}\!\Big((\xi\omega_i)(\eta\omega_i)\Big)\; .
\]
The right hand side is well defined by (\ref{LemI_3_2}).

\rm
By associativity, we may assume $\eta = 1$. We extend the bilinear form $F$-linearly to $F\ts_K L\wr G$ to obtain
\[
\begin{array}{rcl}
(\xi,1)\hspace*{-9mm}
& = & (1\ts\xi,1\ts 1) \\
& = & n^{-1}\,\Tr_F\left((1\ts u)(-)(1\ts \xi)\right) \\
& \aufgl{independence of choice} & n^{-1}\,\Tr_F\left(e_1(-)(1\ts\xi)\right) \\
& \aufgl{algebra isomorphism} & n^{-1}\sum_{i\in [1,k]}\sum_{\theta\in\Theta_i} \Tr_F\Big(((e_1)\omega_{i,\theta})(-)((1\ts\xi)\omega_{i,\theta})\Big) \\
& \aufgl{(\ref{LemI_3_3})} & n^{-1}\sum_{i\in [1,k]}\sum_{\theta\in\Theta_i} 
                             \Tr_F\left((e_1)\omega_{i,\theta}\right)\cdot\Tr_F\left((1\ts\xi)\omega_{i,\theta}\right) \\
& \aufgl{(\ref{LemI_3_2}, \ref{LemI_3_b})} & \sum_{i\in [1,k]}\sum_{\theta\in\Theta_i} (x_i d_i/g)\Tr_F\left(1\liu{\h\theta\,}{\ts}\xi\omega_i\right) \\
& \aufgl{(\ref{LemI_3_c})} & \sum_{i\in [1,k]} (x_i d_i/g)\cdot\tr_{K_i/K}\Tr_{K_i}(\xi\omega_i)\; . \\
\end{array}
\]
\end{Proposition}

\subsection{Conclusions}

\begin{Corollary}
\label{CorI_3_6}
The quotient $x_i d_i/h\in\Z$ is invertible in $K$ for any $i\in [1,k]$. 

\rm
If $x_i d_i/h = 0$ in $K$, then the bilinear form $(-,=)$ on $L\wr G$ would be degenerate by (\ref{PropI_3_5}). This is not the case, and so the result follows. Cf.\ (\ref{LemI_3_2}).
\end{Corollary}

\begin{Corollary}[Fourier inversion]
\label{CorI_3_7}\Absatz
The inverse of the Wedderburn isomorphism is given by
\[
\begin{array}{rcl}
\prodd{i\,\in\, [1,k]} \End_{K_i} X_i & \lraisoa{\omega^{-1}} & L\wr G \\
(\phi_i)_{i\in [1,k]}                 & \lramaps              &
\fracd{1}{n} \sumd{\sigma\,\in\, G,\; l\,\in\, [1,h]} \sigma y_l 
\left[\sumd{i\,\in\, [1,k]}\fracd{x_i d_i}{h}\cdot\tr_{K_i/K}\Tr_{K_i}\left( (y^\ast_l\sigma^{-1})\omega_i\cdot\phi_i\right)\right]\; . \\ 
\end{array}
\]

\rm
Given $\rho\in G$ and $m\in [1,h]$, we verify
\[
\begin{array}{rcl}
(\rho y_m)\omega\omega^{-1}
& = & \fracd{1}{n} \sumd{\sigma\,\in\, G,\; l\,\in\, [1,h]} \sigma y_l 
\left[\sumd{i\,\in\, [1,k]}\fracd{x_i d_i}{h}\cdot\tr_{K_i/K}\Tr_{K_i}\left( (y^\ast_l\sigma^{-1})\omega_i\cdot (\rho y_m)\omega_i\right)\right] \\ 
& \aufgl{(\ref{PropI_3_5})} & \sumd{\sigma\,\in\, G} \sigma \left[\sumd{l\,\in\, [1,h]} y_l (y^\ast_l\sigma^{-1},\rho y_m)\right] \vspace{1mm}\\
& = & \sumd{\sigma\,\in\, G} \sigma \left[\dell_{\rho,\sigma}\sumd{l\,\in\, [1,h]} y_l \Tr_{L/K}(y^\ast_l y_m)\right] \\
& = & \rho y_m\; . \\
\end{array}
\]
\end{Corollary}

Suppose given $i\in [1,k]$, $y\in L$, $\sigma\in G$. With respect to chosen $K_i$-linear bases of $X_i$, we write the endomorphism $(\sigma y)\omega_i$ 
of $X_i$ as a matrix $((\sigma y)\omega_{i;a,b})_{a,b\in [1,x_i]}$.

\begin{Corollary}[Schur relations]
\label{PropI14}\Absatz
Suppose given $i,i'\in [1,k]$, $a,b\in [1,x_i]$, $a',b'\in [1,x_{i'}]$. Then
\[
\frac{g}{x_i d_i}\cdot \dell_{(i;a,b),(i';a',b')} 
\; =\; \sumd{l\in [1,h],\; \sigma\in G} (\sigma y_l)\omega_{i;a,b}\cdot\tr_{K_{i'}/K}\Big( (y_l^\ast\sigma^{-1})\omega_{i';b',a'}\Big)\; .
\]

\rm
By (\ref{CorI_3_7}), we calculate the image of the endomorphism tuple given by the tuple of matrices 
$\Big((\dell_{(i'';a'',b''),(i';a',b')})_{a'',b''\in [1,x_{i''}]}\Big)_{i''\in [1,k]}$ under $\omega^{-1}\omega$ to be
\[
\left(\sumd{l\in [1,h],\;\sigma\in G} (\sigma y_l)\omega_j \frac{x_{i'}d_{i'}}{g} \tr_{K_{i'}/K}\Big((y_l^\ast\sigma^{-1})\omega_{i';b',a'}\Big) \right)_{\!\! j\in [1,k]},
\]
and it remains to compare entries at the position $(i;a,b)$.
\end{Corollary}

\begin{Notation}
\label{NotI14_5}\rm
Given $i\in [1,k]$, we let the {\it character} $\chi_i$ of $X_i$ be defined by
\[
\begin{array}{rcl}
L\wr G & \lraa{\chi_i} &  K \\
\xi    & \lramaps      & \chi_i(\xi)\; :=\; \tr_{K_i/K} \Tr_{K_i}((\xi)\omega_i) \; =\;  \sumd{a\in [1,x_i]}\tr_{K_i/K}((\xi)\omega_{i;a,a})\; . \\
\end{array}
\]
\end{Notation}

\begin{Corollary}[horizontal orthogonality relations]
\label{CorI15}
\Absatz
Suppose given $i,i'\in [1,k]$. Then
\[
g c_i\,\dell_{i,i'} \; =\; \sumd{l\in [1,h],\; \sigma\in G} \chi_i(\sigma y_l) \chi_{i'}(y^\ast_l\sigma^{-1})\; .
\]
\end{Corollary}

\begin{Corollary}[to \ref{CorI_3_7}]
\label{CorI_3_8}\Absatz
The primitive central idempotent $\eps_i$ of $L\wr G$ that belongs to $X_i$, $i\in [1,k]$, is given by 
\[
\eps_i\; =\; \fracd{x_i d_i}{g} \sumd{\sigma\,\in\, G,\; l\,\in\, [1,h]} \sigma y^\ast_l \cdot\chi_i(y_l\sigma^{-1})\; . \\ 
\]
In particular, $\eps_1 = n^{-1}\sumd{\nu\in N}\nu$.

\rm
To evaluate $\eps_1$, it suffices to remark that for $y\in L$ and $\sigma\in G$, the $K$-linear trace of $(y\sigma)\omega_1$ is given by $\{\sigma\in N\}\Tr_{L/K}(y)$, as
to be calculated as $L$-linear trace on $L\ts_K L\iso \prod_{\lambda\in H} L$ using the primitive idempotent basis.
\end{Corollary}

\begin{quote}
\begin{footnotesize}
\begin{Remark}
\label{RemI_3_9}\rm
For the principal module, there is an inversion formula due to {\sc E.\ Noether} {\bf\cite[\rm \S 1]{Noe34}}. Write $\phi\in\End_K L$ as 
$y_l\phi =: \sum_{m\in [1,h]} \phi_{l,m} y_m$ with $\phi_{l,m}\in K$. Then
\[
\begin{array}{rcl}
\End_K L & \lraisoa{\omega_H^{-1}} & L\wr H \hspace*{5mm}(\;\iso\; L\wr G \eps_1 \;\tm\; L\wr G\; ) \\
\phi     & \lramaps                & \sum_{l,m\in [1,h]} \phi_{l,m} y_l^\ast (\sum_{\lambda\in H} \lambda) y_m\; . \\
\end{array}
\]
To see this, we let $\phi_1 = (\rho z)\omega_1$ with $\rho\in G$ and $z\in L$ and apply (\ref{CorI_3_7}) to
\[
\begin{array}{rcl}
(\phi_1,0,\dots,0)
& \lramapsa{\omega^{-1}} & n^{-1}\sum_{\sigma\in G} \sigma\sum_{l\in [1,h]} y_l \Tr_K((y^\ast_l \sigma^{-1}\rho z)\omega_1) \\
& = &  n^{-1}\sum_{\sigma\in G} \sigma\{\sigma^{-1}\rho\in N\}\sum_{l\in [1,h]} y_l \Tr_{L/K}(y^\ast_l z) \\
& = & \rho z\eps_1\; , \\
\end{array}
\]
whereas Noether's formula yields
\[
\begin{array}{rcl}
\phi_1 
& \lramaps & \sum_{l\in [1,h]} y^\ast_l (\sum_{\lambda\in H} \lambda) (y_l\phi_1)  \\
& = & \left(\sum_{\lambda\in H} \lambda \sum_{l\in [1,h]} y^{\ast,\lambda}_l y_l\right) \rho z \; , \\
\end{array}
\]
and it remains to be shown that $\sum_{l\in [1,h]} y^{\ast,\lambda}_l y_l = \dell_{\lambda,1}$. Since the left hand side is independent of the choice of the basis, we may tensor
with $L\ts_K -$ and use the primitive idempotent basis of $L\ts_K L\iso\prod_{\mu\in H} L$, which is self dual. Orthogonality of these idempotents yields the result, for $\lambda\neq 1$
permutes them without fixed points.
\end{Remark}

\begin{Question}
\label{QI_3_10}\rm
Suppose $N = 1$. Given $f\in H^2(G,L^\ast)$, there is a central simple $K$-algebra $L\wr_f G$ (denoted $(L/K,f)$ in {\bf\cite[\rm 29.6]{R75}}). If $f$ is trivial, (\ref{CorI_3_7}) or 
(\ref{RemI_3_9}) allow to deduce an orthogonal primitive idempotent decomposition. I do not know such a decomposition in the general case, nor a workable description of
the simple module (as an $L$-vectorspace with compatible $G$-operation, say). Cf.\ {\bf\cite[\rm 29.22, 32.19]{R75}}, {\bf\cite[\rm 2.4]{W63}}.
\end{Question}
\end{footnotesize}
\end{quote}

\subsection{A Morita equivalence}
\label{SubsecMor}

\begin{Lemma}
\label{LemMor1}
Considering $L\wr G$ as a left module over $LN$, we have an isomorphism
\[
\begin{array}{rcl}
L\ts_K L\wr G & \lraiso  & \End_{LN}(L\wr G) \;\; \iso\;\; (LN)_h \\
y\ts z\sigma  & \lramaps & \Big( w\rho\lramaps y \cdot w\rho\cdot z\sigma \Big)\; \\
\end{array}
\]
of algebras over $L$.

\rm 
Since both source and target are of dimension $gh$ over $L$, it suffices to prove injectivity. Suppose given an element
$\sum_{\sigma\in G,\, i\in I_\sigma} y_{\sigma,i} \ts z_{\sigma,i}\sigma$ that is mapped to zero. From $\sum_{\sigma\in G,\, i\in I_\sigma} y_{\sigma,i} \rho z_{\sigma,i}\sigma = 0$
for all $\rho\in G$ we conclude that $\sum_{i\in I_\sigma} y_{\sigma,i} z_{\sigma,i}^{\rho^{-1}} = 0$ for all $\sigma,\rho\in G$. The isomorphism 
$L\ts_K L\lraiso \prod_{\lambda\in H} L$ now shows that $\sum_{i\in I_\sigma} y_{\sigma,i} \ts z_{\sigma,i} = 0$ for all $\sigma\in G$, as was to be shown.
\end{Lemma}

\begin{Corollary}
\label{CorMor2}
$L\ts_K L\wr G$ is Morita equivalent to $LN$. 

\rm
Concerning the situation over $K$, cf.\ (\ref{CorI17_8}, \ref{ExI2_6}, \ref{ExI2_7}).
\end{Corollary}

\begin{Corollary}
\label{CorMor3}
The field $L(\zeta_{\exp(N)})$ is a splitting field for $L\wr G$ over $K$.

\rm
By (\ref{LemMor1}), we have $L(\zeta_{\exp(N)})\ts_K L\wr G \iso (L(\zeta_{\exp(N)}) N)_h$, and $L(\zeta_{\exp(N)})$ is a splitting field for $N$ by 
{\bf\cite[\rm 12.3, cor.\ to th.\ 24]{Se77}}. Cf.\ also (\ref{ExI2_6}).
\end{Corollary}

\subsection{The center of $L\wr G$}
\label{SubsecCent}

Let $\mb{\rm Cl}^G_N\tm N$ denote a set of representatives of the $G$-orbits of $N$ by conjugation. Given $\sigma\in N$, we let $h_\sigma := h/|C_G(\sigma)N/N| = g/|C_G(\sigma)N|$. Let 
$(y_{\sigma,i})_{i\in [1,h_\sigma]}$ be a $K$-linear basis of $L_\sigma := \mb{Fix}_{C_G(\sigma)}L\tm L$, where $\sigma\in \mb{\rm Cl}^G_N$.

\begin{Lemma}[{\bf\cite[\rm p.\ 113]{NaSh35}}]
\label{LemI16}
The center $Z(L\wr G)$ of $L\wr G$ has the $K$-linear basis
\[
\Big(\sumd{\rho\in C_G(\sigma)\< G} \sigma^\rho y_{\sigma,i}^\rho\Big)_{\sigma\in \mb{\scr\rm Cl}^G_N,\; i\in [1,h_\sigma]}\; .
\]

\rm
Suppose given $\sum_{\sigma\in G}\sigma t_\sigma\in Z(L\wr G)$. For any $y\in L$, we have $y(\sum_{\sigma\in G}\sigma t_\sigma) - (\sum_{\sigma\in G}\sigma t_\sigma)y = 0$, whence 
$\sum_{\sigma\in G}\sigma t_\sigma = \sum_{\sigma\in N}\sigma t_\sigma$. Suppose given $\rho\in G$. From 
\[
\sum_{\sigma\in N}\sigma t_\sigma = (\sum_{\sigma\in N}\sigma t_\sigma)^\rho = \sum_{\sigma\in N}\sigma t_{\!\sigma^{\rho^{-1}}}^\rho,
\]
we conclude that $t_{\sigma^\rho} = t_\sigma^\rho$ for all $\sigma\in N$. Letting $\rho\in C_G(\sigma)$, we conclude that 
$t_\sigma\in L_\sigma$ for all $\sigma\in \mb{\rm Cl}^G_N$.
\end{Lemma}

\begin{Corollary}[{\bf\cite[\rm Satz 2]{NaSh35}}]
\label{CorI17}
The $K$-linear dimension of the center $Z(L\wr G)$, i.e.\ \linebreak[4] $\sum_{l\in [1,k]} \dim_K Z(K_l)$, equals the number of conjugacy classes of $N$.

\rm
{\it First proof.}
We let $G$ act on the set of conjugacy classes of $N$, the $G$-orbits being represented by the $N$-orbits of the elements of $\mb{\rm Cl}^G_N$. The stabilizer of the 
$N$-orbit of $\sigma\in\mb{\rm Cl}^G_N$ in $G$ is given by $C_G(\sigma)N$. But the dimension of $Z(L\wr G)$ over $K$ is given by 
\[
\sumd{\sigma\in\mb{\rm\scr Cl}^G_N} h_\sigma = \sumd{\sigma\in\mb{\rm\scr Cl}^G_N} g/|C_G(\sigma)N|,
\]
which is the sum of the orbit lengths of this operation, indexed by representatives of orbit representatives.

{\it Second proof.} 
We have $\dim_K Z(L\wr G) = \dim_L Z(L\ts_K L\wr G) \aufgl{(\ref{CorMor2})} \dim_L Z(LN)$.
\end{Corollary}

\begin{Corollary}
\label{CorI17_5}
If $K = K_i$ for all $i\in [1,k]$, the number of irreducible modules of $L\wr G$ is given by number of conjugacy classes of $N$.
\end{Corollary}

\begin{Corollary}
\label{CorI17_5_5}
Given $i\in [1,k]$, we have $\chi_i(\sigma y) = 0$ if $\sigma\not\in N$. If $\sigma\in N$, then 
\[
\sum_{l\in [1,h]} y_l^{\ast,\rho} \chi_i(y_l \sigma) \; =\; \sum_{l\in [1,h]} y_l^\ast \chi_i(y_l^{\rho^{-1}} \sigma)\; .
\]
for $\rho\in G$.

\rm
The primitive central idempotent $\eps_i$ described in (\ref{CorI_3_8}) is contained in $Z(L\wr G)$, so that the result follows by (\ref{LemI16}).
\end{Corollary}

\begin{Assumption}
\label{AssI17_6}
For all $\nu\in N$, we have $C_G(\nu)\cdot N = G$.

\rm
This assumption holds if $N = 1$, if $N = G$, or if $N\leq Z(G)$.
\end{Assumption}

\begin{Lemma}[cf.\ {\bf\cite[\rm p.\ 21]{Bo84}}]
\label{LemI17_7}
If (\ref{AssI17_6}) holds, then $Z(KN) = Z(L\wr G)$.

\rm
Consider a basis element $\sum_{\rho\in C_G(\sigma)\< G} \sigma^\rho y_{\sigma,i}^\rho$ of $L\wr G$ as in (\ref{LemI16}).  Since $L_\sigma = K$ by (\ref{AssI17_6}), we may assume
$y_{\sigma,i} = 1$. By (\ref{AssI17_6}), the $G$-orbits and the $N$-orbits on $N$ coincide. Thus our basis element 
equals the class sum $\sum_{\rho\in C_N(\sigma)\< N} \sigma^\rho\in Z(KN)$.
\end{Lemma}

\begin{Corollary}
\label{CorI17_8}
If (\ref{AssI17_6}) holds, and if moreover all endomorphism rings of simple modules of $KN$ and of $L\wr G$ are commutative, then $KN$ and $L\wr G$ are Morita equivalent.

\rm
Concerning the situation after scalar extension over $K$ to $L$, cf.\ (\ref{CorMor2}). Cf.\ also (\ref{ExI2_6}, \ref{ExI2_7}).
\end{Corollary}

\subsection{(Counter)examples}

\begin{quote}
\begin{footnotesize}
\begin{Example}
\label{ExI2_5}\rm %
Suppose given a prime $p\geq 3$, and let $\zeta_{p^2}$ be a primitive $p^2$th root of unity over $\Q$. Let $\pi_{p^2}\in\Q(\zeta_{p^2})$ be defined by 
$\pi_{p^2} := \prod_{j\in [1,p-1]} (\zeta_{p^2}^{j^p} - 1)$. Let $\zeta_p := \zeta_{p^2}^p$.

Consider the operation of $C_{p^2}$ on $\Q(\pi_{p^2})$ given by composition of the surjection $C_{p^2}\lra C_p$ with the operation of the Galois group over $\Q$. Let $\sigma$ be a 
generator of $C_{p^2}$, let $z\in\Q(\pi_{p^2})$. The principal module is given by
\[
\begin{array}{rcl}
\Q(\pi_{p^2})\wr C_{p^2} & \lra     & \End_\ssQ \Q(\pi_{p^2}) \\
\sigma                   & \lramaps & (y\lramaps y^\sigma) \\
z                        & \lramaps & (y\lramaps yz)\; . \\
\end{array}
\]
Moreover, $\Q(\zeta_{p^2})$ becomes a module by means of
\[
\begin{array}{rcl}
\Q(\pi_{p^2})\wr C_{p^2} & \lra     & \End_{\ssQ(\zeta_p)} \Q(\zeta_{p^2}) \\
\sigma                   & \lramaps & (y\lramaps y^\sigma \zeta_{p^2}) \\
z                        & \lramaps & (y\lramaps yz), \\
\end{array}
\]
where $C_{p^2}$ operates on $\Q(\zeta_{p^2})$ via composition of the surjection $C_{p^2}\lra C_p$ with the operation of the Galois group $C_p$ of $\Q(\zeta_{p^2})$ over 
$\Q(\zeta_p)$, so that $\sigma$ maps $\zeta_{p^2}$ to $\zeta_{p^2}^{1+p}$. The element $\sigma^p$ is mapped to the multiplication with
\[
\Nrm_{\ssQ(\zeta_{p^2})/\ssQ(\zeta_p)}(\zeta_{p^2}) \; =\; \zeta_{p^2}^{p+\smatze{p}{2}p} \; =\;  \zeta_{p^2}^p \; .
\]
In particular, $\sigma^{p^2}$ is mapped to the identity.

We claim that $\Q(\zeta_{p^2})$ is a simple module not isomorphic to $\Q(\pi_{p^2})$, more precisely, we claim its endomorphism ring to be equal to $\Q(\zeta_p)$. Since such an endomorphism
is required to commute with $\Q(\pi_{p^2})$ and with $\sigma^p$, it is determined by the image $y$ of $1$. Then commutativity with $\sigma$ yields 
$y = y^\sigma$. As Wedderburn isomorphism, we obtain
\[
\Q(\pi_{p^2})\wr C_{p^2} \;\lraisoa{\omega}\; (\Q)_p\ti (\Q(\zeta_p))_p,
\]
since this surjection is actually an isomorphism for dimension reasons. So in this case $\Q C_p$ and $\Q(\pi_{p^2})\wr C_{p^2}$ are in fact Morita equivalent, as shown
by construction of the progenerator $\Q(\pi_{p^2})\ds \Q(\zeta_{p^2})$ (cf.\ \ref{CorI17_8}).
\end{Example}

\begin{Example}
\label{ExI2_6}\rm
Let $G = \Sl_3$, $N = C_3$, $L = \Q(i)$, $K = \Q$, where $i^{(1,2)} = -i$. We note that (\ref{AssI17_6}) fails, since $C_G((1,2,3)) = N$. According to (\ref{LemI16}), the center
of $\Q(i)\wr\Sl_3$ has the basis $( 1, (1,2,3) + (1,3,2), i(1,2,3) - i(1,3,2))$. Thus we have an isomorphism of $\Q$-algebras
\[
\begin{array}{rclcl}
\Q & \ti & \Q(\sqrt{3}) & \lraiso  & Z(\Q(i)\wr\Sl_3) \\
1  & \ti & 0            & \lramaps & \frac{1}{3}(1 + (1,2,3) + (1,3,2)) \\
0  & \ti & 1            & \lramaps & \frac{1}{3}(2 - (1,2,3) - (1,3,2)) \\
0  & \ti & \sqrt{3}     & \lramaps & i(1,2,3) - i(1,3,2)\; . \\
\end{array}
\]
In particular, $Z(L\wr G)$ is not isomorphic to $Z(KN)$ (cf.\ \ref{LemI17_7} and \ref{CorI17_8}). Moreover, $K(\zeta_{|G|})$ is not a splitting field for the $K$-algebra
$L\wr G$ (but cf.\ \ref{CorMor3}). 
\end{Example}

\begin{Example}
\label{ExI2_7}\rm
Let $G = Q_8 = \spi{i,j\; |\; i^4, j^4, i^2 = j^2, ij = ji^3}$ be the quaternion group with $8$ elements, $N = \spi{i^2}\leq Z(G)$, $L = \Q(\zeta_8)$, $K = \Q$, 
$\zeta_8^i = \zeta_8^{-1}$, $\zeta_8^j = \zeta_8^3$. We remark that (\ref{AssI17_6}) is fulfilled. We have an idempotent
\[
e \; :=\; \frac{1-i^2}{2}\cdot \frac{1 + \zeta_8 j}{2}\; \in\; \Q(\zeta_8)\wr Q_8
\]
and an isomorphism
\[
\begin{array}{rcl}
\HI_\sQ & \lraiso  & e(\Q(\zeta_8)\wr Q_8) e \\
I       & \lramaps & e\zeta_8 i e \\
J       & \lramaps & e (\zeta_8 + i) e\; ,  \\
\end{array}
\]
where $\HI_\sQ = \Q\spi{I,J}/(I^2 + 1, J^2 + 1, IJ + JI)$ is a rational quaternion skewfield. Now, $Z(KN) = Z(L\wr G)$ by (\ref{LemI17_7}), but $KN$ and $L\wr G$ are
not Morita equivalent (cf.\ \ref{CorI17_8}).
\end{Example}
\end{footnotesize}
\end{quote}

\subsection{The central Plancherel formula}
\label{SubSecCB}

We have an isomorphism
\[
\begin{array}{rcl}
Z(L\wr G) & \lraisoa{\omega^Z} & \prodd{i\in [1,k]} Z(K_i) \\
\xi       & \lramaps           & ((\xi)\omega^Z_i)_i := ((\xi)\omega_i)_i\; ,
\end{array}
\]
where $(\xi)\omega_i\in\End_{K_i} X_i$ is actually contained in $Z(K_i)\tm \End_{K_i} X_i$. Recall that $(\eps_j)\omega^Z_i = \dell_{j,i}$ for $i,j\in [1,k]$ (\ref{CorI_3_8}).

Given $\sum_{\rho\in C_G(\nu)\< G} (\nu y)^\rho,\; \sum_{\rho'\in C_G(\nu')\< G} (\nu' y')^{\rho'}\;\in\; Z(L\wr G)$, where $\nu,\nu'\in N$, $y\in L_\nu$ and $y'\in L_{\nu'}\;$, we let
\[
\begin{array}{l}
\Big(\sum_{\rho\in C_G(\nu)\< G} (\nu y)^\rho,\; \sum_{\rho'\in C_G(\nu')\< G} (\nu' y')^{\rho'}\Big)^{\! Z} \vspace{2mm}\\
\; :=\;
\left\{  
\begin{array}{ll}
[N:C_N(\nu)]\cdot\Tr_{L_\nu/K}(yy'^\tau) & \mb{if $\nu = (\nu'^{-1})^\tau$ with $\tau\in G$} \\
0                                        & \mb{if $\nu$ and $\nu'^{-1}$ are not conjugate in $G$} \\
\end{array}
\right. \\
\end{array}
\]
define a $K$-bilinear form on $Z(L\wr G)$ (cf.\ \ref{LemI16}). Note that $yy'^\tau\in L_\nu$. This bilinear form is symmetric and nondegenerate. 

\begin{Remark}
\label{RemCB1}
We have $\; h\cdot (\xi,\xi')^Z\; =\; (\xi,\xi')\;$.

\rm
Let $\nu\in \mb{\rm Cl}^G_N$ and $\nu'\in(\mb{\rm Cl}^G_N)^{-1}$, $y\in L_\nu$ and $y'\in L_{\nu'}\;$. Then
\[
\begin{array}{rcl}
\Big(\sum_{\rho\in C_G(\nu)\< G} (\nu y)^\rho,\; \sum_{\rho'\in C_G(\nu')\< G} (\nu' y')^{\rho'}\Big)
& = & \sum_{\rho\in C_G(\nu)\< G} \dell_{\nu,\nu'^{-1}}\Tr_{L/K}(y^\rho y'^\rho) \\
& = & [G:C_G(\nu)]\cdot\dell_{\nu,\nu'^{-1}}\Tr_{L/K}(yy') \\
& = & [G:C_N(\nu)]\cdot\dell_{\nu,\nu'^{-1}}\Tr_{L_\nu/K}(yy')\; . \\
\end{array}
\]
Thus, if $h$ is invertible in $K$, the following discussion simplifies considerably.
\end{Remark}

\begin{Lemma}
\label{LemCB2}
The bilinear form $(-,=)^Z$ is associative.

\rm
Let $\nu\in \mb{\rm Cl}^G_N$ and $\nu'\in(\mb{\rm Cl}^G_N)^{-1}$, $y\in L_\nu$ and $y'\in L_{\nu'}\;$. Then
\[
\begin{array}{cl}
  & \left(\left(\sum_{\rho\in C_G(\nu)\< G} (\nu y)^\rho\right)\left(\sum_{\rho'\in C_G(\nu')\< G} (\nu' y')^{\rho'}\right),\; 1\right)^{\! Z} \\
= & \dell_{\nu,\nu'^{-1}}\sum_{\rho\in C_G(\nu)\< G} y^\rho y'^\rho \\
= & [N:C_N(\nu)]\cdot \dell_{\nu,\nu'^{-1}}\Tr_{L_\nu/K}(yy')\; . \\
\end{array}
\]
\end{Lemma}

\begin{Proposition}[central Plancherel formula]
\label{PropCB3}
Given $\xi,\xi'\in Z(L\wr G)$, we have
\[
(\xi,\xi')^Z\; =\; \sumd{i\in [1,k]} \frac{1}{n}\left(\frac{x_i d_i}{h}\right)^{\!\! 2}\Tr_{Z(K_i)/K}\left((\xi)\omega^Z_i (\xi')\omega^Z_i\right)\; .
\]

\rm
By associativity, we may assume $\xi' = 1$ (\ref{LemCB2}). Let $\xi =: (\zeta) (\omega^Z)^{-1}$, where 
$\zeta =: (z_i)_{i\in [1,k]}\in\prod_{i\in [1,k]} Z(K_i) \tm \prod_{i\in [1,k]} \End_{K_i} X_i$. Using the auxiliary notation (\ref{NotI_3_1}), and viewing $L\tm F\ts_K L$, we obtain
\[
\begin{array}{cl}
  & ((\zeta) (\omega^Z)^{-1},1)^Z \\
\aufgl{(\ref{CorI_3_7})} & \left(\sum_{\sigma\in G,\; l\in [1,h]} \sigma y_l
                           \left[\sum_{i\in [1,k]}\frac{x_i d_i}{g}\cdot\tr_{K_i/K}\Tr_{K_i}((y_l^\ast\sigma^{-1})\omega_i\cdot z_i)\right],\; 1 \right)^{\!\! Z} \\
\aufgl{definition}       & \sum_{i\in [1,k]}\sum_{l\in [1,h]} y_l\cdot\frac{x_i d_i}{g}\cdot\tr_{K_i/K}\Tr_{K_i}((y_l^\ast)\omega_i\cdot z_i) \\
\aufgl{(\ref{LemI_3_c}, \ref{LemI_3_b})} & \sum_{i\in [1,k]}\sum_{l\in [1,h]}\sum_{\theta\in\Theta_i} (1\ts y_l)\cdot\frac{x_i d_i}{g}\cdot
                           \Tr_F((1\ts y_l^\ast)\omega_{i,\theta}\cdot (1\liu{\h\theta\,}{\ts} z_i)) \\
\aufgl{\shortstack[c]{independence\\of choice}} & \sum_{i\in [1,k]}\sum_{\mu\in H}\sum_{\theta\in\Theta_i} e_\mu\cdot\frac{x_i d_i}{g}\cdot
                           \Tr_F((e_\mu)\omega_{i,\theta}\cdot (z_i\theta\liu{\h\theta\,}{\ts} 1)) \\
\aufgl{(\ref{LemI_3_2})} & \sum_{i\in [1,k]}\sum_{\mu\in H}\sum_{\theta\in\Theta_i} e_\mu\cdot\frac{x_i d_i}{g}\cdot\frac{x_i d_i}{h}\cdot z_i\theta \\
= & \sum_{i\in [1,k]} \frac{1}{n}\left(\frac{x_i d_i}{h}\right)^{\! 2}\Tr_{Z(K_i)/K}(z_i)\; . \\
\end{array}
\]
\end{Proposition}
\section{Locally integral}
\label{SubsecResInt}

\begin{Notation}
\label{SetupI3_05}\rm
Let $S\tm K$ be a discrete valuation ring with field of fractions $K$, maximal ideal generated by $\pi$ and residue field $\b S := S/\pi S$ of characteristic $p \geq 0$. The $\pi$-adic 
valuation of an element $y\in K$ is denoted by $v(y)$. Let $T\tm L$ be the integral closure of $S$ in $L$. Note that $T$ is a principal ideal
domain. Let $\pi^s S := \Tr_{L/K}(T)$. We have $\pi T = \prod_{i\in [1,d]} \qfk_i^e$, where the $\qfk_i\tm T$ denote the maximal ideals, forming a single Galois orbit, and where 
$e$ is the {\it ramification index} of $\pi$ in $T$ {\bf\cite[\rm I.\S 9]{N92}}. Recall that $\pi$ is said to be {\it tamely ramified} in $T$ if $T/\qfk_i$ is separable over $\b S$ 
for some (hence for all) $i\in [1,d]$, and if, in addition, the ramification index $e$ of $\pi$ in $T$ is not divisible by $p$. Otherwise, it is said to be {\it wildly ramified.}
\end{Notation}

\subsection{Projectivity}

\begin{Lemma}[Maschke locally]
\label{LemI3_1}
For $S$-free $T\wr G$-modules $X$ and $Y$, we have \linebreak[4] \mb{$\pi^s n\cdot\Ext_{T\wr G}^1(Y,X) = 0$.} In particular, if $\eps$ is a central idempotent of $L\wr G$, we have 
$\pi^s n \eps \in T\wr G$.

\rm
Let $u_0\in T$ be such that $\Tr_{L/K}(u_0) = \pi^s$. Consider an extension 
\[
0\;\lra\; X\;\lra\; E\;\lraa{f}\; Y\;\lra\; 0 
\]
of $T\wr G$-modules. We choose a $T$-linear coretraction $Y\lraa{i_0} E$ to $E\lraa{f} Y$ and let 
\[
\begin{array}{rcl}
Y & \lraa{i} & E \\
y & \lramaps & \sumd{\sigma\in G} (y\sigma)i_0\cdot u_0\sigma^{-1}\; , \\
\end{array}
\]
which is $T\wr G$-linear and satisfies $if = \pi^s n\cdot 1_Y$.

Applying this construction to the epimorphism $T\wr G\lraa{f} \eps T\wr G$ given by left multiplication with $\eps$, the map $i$ sends $\eps$ to an element of the 
form $\xi = \xi\eps\in T\wr G$. Composition with $f$ yields $\pi^s n\eps = \eps\xi\eps \aufgl{$\eps$ central} \xi\eps\in T\wr G$.
\end{Lemma}

\begin{Corollary}
\label{CorI3_1_5}
Given $i\in [1,k]$, $\sigma\in G$ and $l\in [1,h]$, we have $\pi^s\frac{x_i d_i}{h} y^\ast_l\chi_i(y_l\sigma)\in T$.

\rm
Apply (\ref{LemI3_1}) to (\ref{CorI_3_8}).
\end{Corollary}

\begin{Proposition} [{{\sc Speiser,} {\sc E.\ Noether,} see {\bf\cite[\rm I.\S 3, Th.\ 3]{F83}}}]
\label{RemI3_2}
\Absatz
The following assertions are equivalent.
\begin{itemize}
\item[(i)] As a module over $T\wr G$, $T$ is projective.
\item[(ii)] As a module over $SG$, $T$ is projective.
\item[(iii)] The prime $p$ does not divide $n$, and the trace $\Tr_{L/K}$ maps $T$ onto $S$.
\item[(iv)] The prime $p$ does not divide $n$, and $\pi$ is tamely ramified in $T$. 
\end{itemize}

\rm
Assertion (i) implies (ii) since $T\wr G$ is free over $SG$. We claim that (ii) implies (i). Given a split $SG$-linear epimorphism $(SG)^{(i)}\lra T$, we 
obtain a split $T\wr G$-linear epimorphism $(T\wr G)^{(i)}\lra T\ts_{SG} T\wr G$ when tensoring with $T\wr G$ over $SG$. But the 
$T\wr G$-linear epimorphism (the counit of the tensor adjunction)
\[
\begin{array}{rclcl}
T & \ts_{SG} & T\wr G   & \lra     & T \\
y & \ts      & \sigma z & \lramaps & y^\sigma z \\
1 & \ts      & z        & \llamaps & z \\
\end{array}
\]
is $T\wr G$-linearly split as indicated.

We claim equivalence of (i) and (iii). As a module over $T\wr G$, $T$ is projective if and only if the $T\wr G$-linear epimorphism 
\[
\begin{array}{rcl}
T\wr G   & \lraa{\eta} & T \\
1        & \lramaps    & 1 \\
\end{array}
\]
is split by some coretraction $T\wr G\llaa{\phi} T$. Amongst the $T$-linear maps from $T$ to $T\wr G$, those are $T\wr G$-linear that send 
$1$ to $y\sum_{\sigma\in G}\sigma$ for some $y\in T$. But this map is a coretraction to $\eta$ if and only if $n\cdot\Tr_{L/K} y = 1$.

We claim equivalence of (iii) and (iv). Reducing modulo $\pi$, the surjectivity of the trace in question is
equivalent to the nonvanishing of the $\b S$-linear trace on some factor $T/\qfk_i^e$ of $\b T$. This trace in turn equals $e$ times the
trace on $T/\qfk_i$, by a filtration argument. Finally, the trace on $T/\qfk_i$ does not vanish if and only if $T/\qfk_i$ is separable over $\b S$.
\end{Proposition}

\begin{Remark}
\label{RemI_3_2_5}\Absit\rm
\begin{itemize}
\item[(i)] If the equivalent conditions of (\ref{RemI3_2}) are fulfilled, then $T\iso_{SG} SH$ {\bf\cite[\rm 32.1]{CR81}}.
\item[(ii)] If $N = 1$, another proof of (\ref{RemI3_2}), (iii $\Rightarrow$ i) can be obtained by {\bf\cite[\rm 28.7]{CR81}}, and by noting that 
$T\iso (\sum_{\sigma\in G} \sigma) T\wr G\tm T\wr G$ is a right ideal.
\end{itemize}
\end{Remark}

\subsection{The Wedderburn colength}

\begin{quote}
\begin{footnotesize}
We give a formula for the colength of the Wedderburn embedding of twisted group rings (\ref{PropI7}), following basically the method of {\sc Plesken} {\bf\cite{P83}} (that 
{\sc G.\ Nebe} taught me). 
\end{footnotesize}
\end{quote}

\begin{Notation}
\label{SetupI3_3}\rm
For each $i\in [1,k]$, we choose a finitely generated $T\wr G$-module $V_i\tm X_i$ such that $KV_i = X_i$ and such that $S_i := \End_{T\wr G} V_i$ is a maximal order in $K_i$
(cf.\ \ref{RemI4}). In particular, we choose $V_1 = T \tm L = X_1$, so that $S_1 = S$. Note that $S_i$ is finitely generated free as a module over $S$, and that $KS_i = K_i$, since 
$S_i = \End_S V_i\cap K_i\tm\End_K X_i$. Note that $V_i$ is free as a module over $S_i$ {\bf\cite[\rm 18.7 (i)]{R75}}, of rank $x_i$. We fix $S_i$-linear bases for $V_i$, $i\in [1,k]$, 
and write $S_i$- (resp.\ $K_i$-) linear endomorphisms of $V_i$ (resp.\ of $X_i$) as matrices over $S_i$ (resp.\ over $K_i$), i.e.\ as elements of $(S_i)_{x_i}$ 
(resp.\ of $(K_i)_{x_i}$). In particular, we identify $\End_{S_i} V_i = (S_i)_{x_i}$.
\end{Notation}

\begin{Remark}
\label{RemI4}\rm
The required choice of $V_i\tm X_i$ can be achieved as follows. Let $S_i$ be a maximal order in $K_i$ (\ref{LemI2_5}), {\bf\cite[\rm 7.18 (ii), 10.4]{R75}}. The $S$-subalgebra $\Xi$ of 
$\End_K X_i$ generated by $S_i$ and the image of $T\wr G$ therein, is finitely generated as a module over $S$. Let $V_i$ be a finitely generated $\Xi$-submodule of $X_i$
such that $KV_i = X_i$. Since $S_i\tm\End_{T\wr G} V_i\tm K_i$, maximality of $S_i$ yields $S_i = \End_{T\wr G} V_i$.
\end{Remark}

\begin{Notation}
\label{NotI_5}\rm
Given a finitely generated torsion $S$-module $M$, we denote by $l(M)$ the length of $M$ in the sense of Jordan-H\"older. 
Given an inclusion of $S$-modules $N\tm M$ such that $M/N$ is a finitely generated torsion $S$-module, the {\it colength of $N$ in $M$} is defined to be the length 
$l(M/N)$, i.e.\ the multiplicity of $\b S$ in $M/N$.

We write shorthand 
\[
\begin{array}{rcl}
A       & := & L\wr G \\
\Lambda & := & T\wr G\; . \\
\end{array}
\]
The $K$-bilinear form on $A$ defined in (\ref{DefI2_6}) restricts to an $S$-bilinear form on $\Lambda$. For an $S$-submodule $U\tm A$, we let
\[
U^\sharp \; :=\;  \{ y\in A\;|\; (y,U)\tm S \}\;\tm\; A
\]
be its {\it dual} $S$-submodule. 
\end{Notation}

\begin{Lemma}
\label{LemI4}
For a finitely generated $S$-submodule $U$ of $A$ such that $KU = A$, the module $U^\sharp$ is also a finitely generated $S$-submodule of $A$ such that $KU^\sharp = A$. 
In this case, we obtain
\[
U\; =\; U^{\sharp\sharp}\; .
\]
Therefore $U\lramaps U^\sharp$ is an antiinvolution of the lattice $\Ul^\mb{\scr f}$ of finitely generated $S$-submodules of $A$ such that $KU = A$.
In particular, given $U,V\in\Ul^\mb{\scr f}$, $U\tm V$, we have $l(V/U) = l(U^\sharp/V^\sharp)$.

\rm
We choose two $S$-linear bases $(u_i)_{i\in [1,gh]}$ and $(u'_i)_{i\in [1,gh]}$ of $U$ such that $(u_i,u_j') = \dell_{i,j} z_i$ for 
some $z_i\in K^\ast$, $i,j\in [1,gh]$. We obtain $S$-linear bases $(z_i^{-1}u_i)_{i\in [1,gh]}$ and $(z_i^{-1}u'_i)_{i\in [1,gh]}$ of $U^\sharp$. With respect to 
these bases, we have $(z_i^{-1} u_i, z_j^{-1} u_j') = \dell_{i,j} z_i^{-1}$, $i,j\in [1,gh]$. Hence by the same argument, applied to this situation, we obtain 
$U^{\sharp\sharp} = U$.
\end{Lemma}

\begin{Notation}
\label{NotI4_3}\rm
Suppose given $i\in [1,k]$. Let 
\[
S_i^+ \; :=\; \{ y\in K_i\; |\; \tr_K(yS_i)\tm S\}\; .
\]
The colength of $S_i$ in $S_i^+$ is denoted by
\[
\delr_{S_i/S} \; :=\; l(S_i^+/S_i)\; .
\]
Furthermore, let 
\[
\begin{array}{rcl}
T^+              & := & \{y\in L\; |\; \Tr_{L/K}(yT)\tm S\} \\
\delr_{T/S}      & := & l(T^+/T) \\
\Lambda^+        & := & \{\sumd{\sigma\in G} \sigma y_\sigma \in A\; |\; y_\sigma\in T^+\}\; . \\
\end{array}
\]
Cf.\ {\bf\cite[\rm III.2.1]{N92}}.
\end{Notation}

\begin{Remark}
\label{LemI4_5}
We have
\[
l(\Lambda^+/\Lambda) \; =\; g\delr_{T/S}\; .
\]
\end{Remark}

\begin{Lemma}
\label{LemI5}
We have
\[
\Lambda^\sharp \; =\;  \Lambda^+\; .
\]

\rm
Given $\sumd{\sigma\in G} \sigma y_\sigma \in A$, $y_\sigma\in L$, we obtain  
\[
(\sumd{\sigma\in G} \sigma y_\sigma, \rho^{-1}z) \; =\; \Tr_{L/K}(y_\rho^{\rho^{-1}} z)\; ,
\]
which is in $S$ for all $\rho\in G$ and for all $z\in T$ if and only if $\sum_{\sigma\in G} \sigma y_\sigma$ is in $\Lambda^\sharp$, but also if and only if it is in $\Lambda^+$.
\end{Lemma}

\begin{Notation}
\label{NotI_5_5}\rm
Let 
\[
\begin{array}{rclcl}
\Gamma   & := & \left(\prodd{i\in [1,k]} (S_i)_{x_i}\right)\omega^{-1}            & \tm & A \\
\Gamma^+ & := & \left(\prodd{i\in [1,k]} \fracd{h}{x_i d_i} (S_i^+)_{x_i}\right)\omega^{-1} & \tm & A\; , \\
\end{array}
\]
where $(S_i^+)_{x_i}$ is the $S$-submodule of $(K_i)_{x_i}$ consisting of matrices entrywise in $S_i^+$. Cf.\ (\ref{LemI_3_2}, \ref{CorI_3_6}).
\end{Notation}

\begin{Lemma}
\label{LemI6}
We have
\[
\Gamma^\sharp \; =\; \left(\prodd{i\in [1,k]} \frac{g}{x_i d_i} (S_i^+)_{x_i}\right)\omega^{-1}\; .
\]
Hence 
\[
\begin{array}{rcl}
l(\Gamma^+/\Gamma)        & = & \sumd{i\in [1,k]} x_i^2\left(\delr_{S_i/S} + v(x_i d_i/h) r_i \right) \\
l(\Gamma^+/\Gamma^\sharp) & = & v(n)gh\; . \\
\end{array}
\]

\rm
We use (\ref{PropI_3_5}) to apply the bilinear form on the right hand side of the Wedderburn isomorphism to a tuple of matrices with only one nonzero matrix which in turn has only 
one nonzero entry, and an arbitrary tuple.
\end{Lemma}

\begin{Theorem}
\label{PropI7}
The colength of the Wedderburn embedding
\[
T\wr G \;\;\hraa{\omega}\; \prodd{i\in [1,k]} \End_{S_i} V_i
\]
is given by
\[
\frac{1}{2}\left(g(\delr_{T/S} + v(n)h) - \sumd{i\in [1,k]} x_i^2 (\delr_{S_i/S} + v(x_i d_i/h) r_i ) \right).
\]
Concerning notation, cf.\ (\ref{SetupI0_5}, \ref{SetupI2_7}, \ref{SetupI3_05}, \ref{SetupI3_3}, \ref{NotI_5}, \ref{NotI_5_5}).

\rm
Using (\ref{LemI5}), the diagram 
\begin{center}
\begin{picture}(550,450)
\put(   0, 200){$\Gamma^\sharp$}
\put(  60, 210){\vector(1,0){120}}
\put(  60, 220){\oval(20,20)[l]}
\put( 200, 200){$\Lambda^\sharp$}
\put( 250, 200){$=$}
\put( 300, 200){$\Lambda^+$}
\put( 300, 400){$\Lambda$}
\put( 310, 370){\vector(0,-1){120}}
\put( 320, 370){\oval(20,20)[t]}
\put( 360, 410){\vector(1,0){120}}
\put( 360, 420){\oval(20,20)[l]}
\put( 500, 400){$\Gamma$}
\put( 510, 370){\vector(0,-1){320}}
\put( 520, 370){\oval(20,20)[t]}
\put( 500,   0){$\Gamma^+$}
\put(  60, 180){\vector(3,-1){420}}
\put(  60, 190){\oval(15,20)[l]}
\end{picture}
\end{center}
of $S$-submodules of $A$ shows that
\[
l(\Lambda^\sharp/\Gamma^\sharp) - l(\Lambda^+/\Lambda) + l(\Gamma/\Lambda) 
+ l(\Gamma^+/\Gamma) - l(\Gamma^+/\Gamma^\sharp) \; =\; 0\; .
\]
The result follows by (\ref{LemI4}, \ref{LemI4_5}, \ref{LemI6}).
\end{Theorem}

\begin{quote}
\begin{footnotesize}
\begin{Question}
\label{QuI7_5}\rm
Given $i\in [1,k]$, it would be desirable to know the colength of the embedding of the quasiblock $T\wr G\eps_i$ into $\End_{S_i} V_i$ via $\omega_i$.
\end{Question}
\end{footnotesize}
\end{quote}

\begin{Corollary}
\label{CorI8}
If $G$ acts faithfully on $L$, the colength of the Wedderburn embedding $T\wr G\hraa{\omega} \End_S T$ is given by
\[
\frac{1}{2}\, g\, \delr_{T/S}\; .
\]
\end{Corollary}

\begin{Corollary}[cf.\ {\bf\cite[\rm 1.1.5]{Kue99}}]
\label{CorI9}
If $G$ acts trivially on $L$, the colength of the Wedderburn embedding $SG\;\hraa{\omega}\prodd{i\in [1,k]}\End_{S_i} V_i$ is given by
\[
\frac{1}{2}\left(g v(g) - \sumd{i\in [1,k]} x_i^2 (\delr_{S_i/S} + v(x_i d_i) r_i ) \right)\; .
\]
\end{Corollary}

\subsection{Some estimates}

\begin{Notation}
\label{Setup9_5}\rm
Suppose the basis $(y_i)_{i\in [1,h]}$ chosen in (\ref{SetupI0_5}) to be an $S$-linear basis of $T$. Let
\[
S\pi^t\; :=\; \{ z\in S\; |\; z y_i^\ast\in T\;\mb{ for all $i,j\in [1,h]$}\}\; .
\]
So $t$ is the maximal elementary divisor exponent of the Gram matrix of the trace bilinear form on $T$ induced by $\Tr_{L/K}$. In particular, $s\leq t$
\end{Notation}

\begin{Lemma}
\label{CorI11}
The cokernel of the Wedderburn embedding
\[
T\wr G \;\;\hraa{\omega} \prodd{i\in [1,k]}\;\End_{S_i} V_i
\]
is annihilated by $n\pi^t$.

\rm
This follows by Fourier inversion (\ref{CorI_3_7}) and by (\ref{LemI_3_2}).
\end{Lemma}

\begin{Lemma}
\label{LemI15_1}
Given $i\in [1,k]$, we have $v(x_i) + v(d_i) \leq v(g) + t$. Cf.\ (\ref{CorI_3_6}).

\rm
The multiple 
\[
\frac{\pi^t g}{x_i d_i}\cdot \eps_i \; =\; \pi^t \sumd{\sigma\,\in\, G,\; l\,\in\, [1,h]} \sigma y^\ast_l \cdot\chi_i(y_l\sigma^{-1})\;\in\; T\wr G
\]
of the central primitive idempotent $\eps_i$ is mapped under $\omega_i$ to $(\pi^t g)/(x_i d_i)$ times the identity on $V_i$ (\ref{CorI_3_8}).
\end{Lemma}

\subsection{The central colength}

\begin{Notation}
\label{NotZ1}\rm
Given $\nu\in N$, we let $T_\nu := \mb{Fix}_{C_G(\nu)} T = T\cap L_\nu$. Note that $T_\nu$ is the integral closure of $S$ in $L_\nu$. Let 
$T_\nu^+ := \{ y\in L_\nu\; |\; \Tr_{L_\nu/K}(y T_\nu)\tm S\}$, let $\delr_{T_\nu/S} := l(T_\nu^+/T_\nu)$. 
Let $Z(S_i)^+ := \{ z\in Z(K_i)\; |\; \Tr_{Z(K_i)/K}(z Z(S_i))\tm S\}$, let $\delr_{Z(S_i)/S} := l(Z(S_i)^+/Z(S_i))$. 
\end{Notation}

\begin{Lemma}
\label{LemZ2}
With respect to the bilinear form $(-,=)^Z$ on $Z(A)$ (section \ref{SubSecCB}), we obtain 
\[
\begin{array}{rcl}
Z(\Lambda)^\sharp
& := & \{\zeta\in Z(A)\; |\; (\zeta,Z(\Lambda))^Z\tm S\} \\
& =  & \left\{\sumd{\sigma\in\mb{\scr\rm Cl}^G_N}\;\sumd{\rho\in C_G(\sigma)\< G} (\sigma y_\sigma)^\rho\;\left|\rule[7mm]{0mm}{0mm}\right.\; 
              y_\sigma\in [N:C_N(\sigma)]^{-1} T^+_\sigma\right\}\; . \\
\end{array}
\]
In particular, the $S$-linear colength of $Z(\Lambda)$ in $Z(\Lambda)^\natural$ is given by
\[
\sum_{\nu\in\mb{\scr\rm Cl}^G_N} \delr_{T_\nu/S} + v([N:C_N(\nu)])h_\nu\; .
\]

\rm
Suppose given $\zeta = \sum_{\sigma\in\mb{\scr\rm Cl}^G_N}\sum_{\rho\in C_G(\sigma)\< G} (\sigma y_\sigma)^\rho\in Z(A)$, where $y_\sigma\in L_\sigma$ (\ref{LemI16}), and suppose given
$\sum_{\rho'\in C_G(\nu)\< G} (\nu^{-1} z)^{\rho'}\in Z(\Lambda)$, where $\nu\in\mb{\rm Cl}^G_N$ and $z\in T_\nu$. We obtain
\[
\left(\sum_{\sigma\in\mb{\scr\rm Cl}^G_N}\sum_{\rho\in C_G(\sigma)\< G} (\sigma y_\sigma)^\rho\; , \sum_{\rho'\in C_G(\nu)\< G} (\nu^{-1} z)^{\rho'}\right)^{\!\! Z} 
\; =\; [N : C_N(\nu)]\cdot\Tr_{L_\nu/K}(y_\nu z)\; .
\]
\end{Lemma}

\begin{Lemma}
\label{LemZ3}
With respect to the bilinear form $(-,=)^Z$ on $Z(A)$ (cf.\ section \ref{SubSecCB}), we obtain
\[
\begin{array}{rcl}
Z(\Gamma)^\sharp
& := & \{\zeta\in Z(A)\; |\; (\zeta,Z(\Gamma))^Z\tm S\} \\
& =  & \{((z_i)_{i\in [1,k]})(\omega^Z)^{-1}\; |\; z_i\in \frac{gh}{x_i^2 d_i^2}Z(S_i)^+\}\; . \\
\end{array}
\]

\rm
This follows by (\ref{PropCB3}).
\end{Lemma}

\begin{Proposition}
\label{PropZ4}
The colength of the central Wedderburn embedding 
\[
Z(T\wr G)\;\hraa{\omega^Z}\; \prodd{i\in [1,k]} Z(S_i)
\]
is given by
\[
\frac{1}{2}\left(\left(\sum_{\nu\in\mb{\scr\rm Cl}^G_N} \delr_{T_\nu/S} + v([N : C_N(\nu)])h_\nu\right) 
- \left(\sum_{i\in [1,k]}\delr_{Z(S_i)/S} + v\!\left(\frac{x_i^2 d_i^2}{gh}\right) c_i\right)\right)\; .
\]

\rm
This follows by (\ref{LemZ2}, \ref{LemZ3}) by the argument of (\ref{PropI7}).
\end{Proposition}

\begin{Corollary}[{cf.\ {\bf\cite[\rm 4.1]{CPW87}}}]
\label{CorZ5}
If $G$ acts trivially on $L$, the colength of the central Wedderburn embedding 
$\;
Z(T\wr G)\;\hraa{\omega^Z}\;\prodd{i\in [1,k]} Z(S_i)
\;$
is given by
\[
\frac{1}{2}\left(\left(\sum_{\sigma\in\mb{\scr\rm Cl}^G_G} v([G : C_G(\sigma)])\right) 
- \left(\sum_{i\in [1,k]}\delr_{Z(S_i)/S} + v\!\left(\frac{x_i^2 d_i^2}{g}\right) c_i\right)\right)\; .
\]
\end{Corollary}
\section{Modular}

\subsection{Brauer-Nesbitt}

\begin{Corollary}[to \ref{PropI14}]
\label{CorBN1}
Let $i\in [1,k]$. Suppose that $K = K_i$ and that $v(x_i) = v(g) + t$ (cf.\ \ref{Setup9_5}, \ref{LemI15_1}). Then the reduction $V_i/\pi V_i$ 
is a simple $\b T\wr G$-module.

\rm
Assume $V_i/\pi V_i$ to have a nontrivial submodule. Choose an $S$-linear basis of $V_i$ such that $v\left( (\sigma y_l)\omega_{i;\, x_i,1}\right) > 0$ for all $\sigma\in G$, 
$l\in [1,h]$. Then
\[
\frac{g}{x_i} \;\aufgl{(\ref{PropI14})}\; \sumd{l\in [1,h],\;\sigma\in G} (\sigma y_l)\omega_{i;\, x_i,1}\cdot (y^\ast_l\sigma^{-1})\omega_{i;\, 1,x_i}\; .
\]
Now the left hand side has valuation $-t$, whereas the right hand side has valuation $> -t$. This contradicts our assumption, and therefore $V_i/\pi V_i$ is simple.
\end{Corollary}

\begin{Remark}
\label{RemBN2}
The $\b S$-algebra $\b T\wr G$ is semisimple if and only if $e = 1$, $T/\qfk_1$ is separable over $\b S$ and $p$ does not divide $n$. We have $e = 1$ if and only if $t = 0$.

\rm
For the necessity of these three conditions, we consider the split epimorphisms $\b T\lra \b T_0$, which implies that $e = 1$, and $\b T\wr G\lra \b T$, which then implies that 
$T/\qfk_1$ is separable over $\b S$ and that $p$ does not divide $n$ (cf.\ \ref{RemI3_2}).

First of all, if $e = 1$ and $T/\qfk_1$ is separable over $\b S$, then $\Tr_{\b T/\b S}(\b T) = \b S$ (cf.\ \ref{RemI3_2}).
Sufficiency now follows as in (\ref{LemI3_1}), using an element $\b u_0\in \b T$ such that $\Tr_{\b T/\b S}(\b u_0) = 1$.

We have $t = 0$ if and only if $T^+ = T$, i.e.\ if and only if the discriminant of $T/S$ is invertible in $S$. By {\bf\cite[\rm III.2.12]{N92}}, this holds if and only if $e = 1$.
\end{Remark}

\begin{Remark}
\label{RemBN3}
Let $i\in [1,k]$. If $\b T\wr G$ is semisimple and if $K = K_i$, then the reduction $V_i/\pi V_i$ is a simple $\b T\wr G$-module.

\rm
By (\ref{CorBN1}), we have to show that $v(x_i) = v(h) + v(n) + t$. Now $t = 0$ and $v(n) = 0$ by (\ref{RemBN2}), $v(x_i) \geq v(h)$ by (\ref{LemI_3_2}) and $v(x_i) \leq v(h)$
by (\ref{LemI15_1}).
\end{Remark}

\begin{quote}
\begin{footnotesize}
\begin{Question}
\label{QBN4}\rm
Suppose $\b T_0\wr G$ to be semisimple. For instance, this is the case if $\b T\wr G$ is semisimple (\ref{RemBN2}). Then $\Jac(\b T\wr G) = \Jac(T)(\b T\wr G)$, being a nilpotent ideal
with semisimple quotient. Are the isoclasses of simple $\b T\wr G$-modules represented by $V_i/\Jac(T) V_i$ for \mb{$i\in [1,k]$}~? Cf.\ {\bf\cite[\rm \S 28, ex.\ 5]{CR81}}.
\end{Question}
\end{footnotesize}
\end{quote}

\subsection{Counting simple modules}

\begin{quote}
\begin{footnotesize}
The arguments in this section are an adaption of {\sc Brauer'}s proof of the classical case \mb{\bf\cite[\rm 3B]{B56}}. 
\end{footnotesize}
\end{quote}

\begin{Notation}
\label{NotMod1}\rm
Recall from (\ref{SetupI3_05}) that we write $\b S = S/\pi S$.  

\fbox{Assume $\b S$ to be finite, and let $\b S =: p^a = : q$.} 

On the other hand, in this section we may admit $n\cdot 1_K = 0$. 

Denote $\b T := T/\pi T$ and $\b\qfk_i := \qfk_i/\pi T$. Note that 
$\b T\iso\prod_{i\in [1,d]} \b T_{\qfk_i}\iso\prod_{i\in [1,d]} T/\qfk_i^e$. 

Let $\bL := \b T\wr G$, so $A\om\Lambda\ra\b\Lambda$. Let $\bLC$ be the $\b S$-subspace of $\bL$ generated by the elements $ab - ba$, for $a,b\in\bL$. Let 
\[
\bLCP \; :=\; \{ \xi\in\bL\; |\; \xi^{p^m}\in\bLC\;\mb{ for some $m\geq 0$}\}\; .
\]
\end{Notation}

\begin{Lemma}[\bf\cite{Kuel81}]
\label{LemOmni}\Absit
\begin{itemize}
\item[(i)] Given $\xi,\eta\in\bL$ and $m\geq 0$, we have $(\xi + \eta)^{p^m}\con_\bLC \xi^{p^m} + \eta^{p^m}$. 
\item[(ii)] The set $\bLCP$ is an $\b S$-linear subspace of $\bL$. There exists $M\geq 0$ such that for any $m\geq M$
  \[
  \bLCP \; =\; \{ \xi\in\bL\; |\; \xi^{p^m}\in\bLC\}\; =: \bLC^{p^{-m}}\; .
  \] 
\item[(iii)] Let $\mu\geq 0$ be such that $a\mu\geq M$ and such that $a\mu\geq v_p(o(\sigma))$ for all $\sigma\in G$. We have an injective $\b S$-linear map
  \[
  \begin{array}{rcl}
  \bL/\bLCP\; =\; \bL/\bLC^{q^{-\mu}} & \hraa{\iota} & \bL/\bLC \\
                               \xi    & \lramaps     & \xi^{q^\mu}\; . \\
  \end{array}
  \] 
\item[(iv)] We have
  \[
  \dim_{\b S} Z(\bL/\Jac(\bL)) \; =\; \dim_{\b S}(\bL/\bLCP)\; .
  \] 
\end{itemize}
\end{Lemma}

\begin{Notation}
\label{NotMod7}\rm
Let $\b T_0 := \b T/\Jac(\b T) = T/\Jac(T)\iso\prod_{i\in [1,d]} \b T_{0,\qfk_i}\iso\prod_{i\in [1,d]} T/\qfk_i$. 
Deviating from former notation, we let $\ClG$ be a system of representatives of conjugacy classes of $G$, and let $\ClGpp := \{\sigma\in\ClG \; |\; o(\sigma)\not\con_p 0\}\tm\ClG$ 
be the subset of $p'$-representatives.

Given $\sigma\in\ClG$, we define the subspace
\[
\fbox{$\;\; 
V_\sigma \; :=\; \b S\spi{(y^\sigma - y)z,\; y^\rho - y\; |\; y,z\in\b T_0,\; \rho\in C_G(\sigma)}\; \tm\; \b T_0\; .
\;\;$}
\]
\end{Notation}

\begin{Lemma}
\label{LemMod8}
There is an $\b S$-linear map
\[
\begin{array}{rcl}
\bL/\bLC         & \lraa{\phi} & {\dis\Ds_{\kappa\in\sClG}} \b T_0/V_\kappa \\
\sigma^\rho y    & \lramaps    & (\dell_{\sigma,\kappa} y^{\rho^{-1}})_\kappa\; , \\
\end{array}
\]
where $y\in\b T$, $\sigma\in\ClG$ and $\rho\in G$.

\rm
We need to show the independence of the choice of $\rho$. But if $\tau\in C_G(\sigma)$, then \linebreak[4] $y^{\rho^{-1}} - y^{\rho^{-1}\tau^{-1}}\in V_\sigma$.

We need to show that given $\tau,\w\tau\in G$, $y,\w y\in\b T$, the element $\tau y\w\tau\w y - \w\tau\w y\tau y$ is mapped to zero. Writing $\tau\w\tau = \sigma^\rho$
with $\sigma\in\ClG$ and $\rho\in G$, we get
\[
\begin{array}{rcl}
(\tau y\w\tau\w y - \w\tau\w y\tau y)\phi
& = & (\sigma^\rho y^{\w\tau} \w y - \sigma^{\rho\tau} \w y^\tau y)\phi \\
& = & (\dell_{\sigma,\kappa} y^{\w\tau\rho^{-1}} \w y^{\rho^{-1}} - \w y^{\rho^{-1}} y^{\tau^{-1}\rho^{-1}})_\kappa \\
& = & (\dell_{\sigma,\kappa} (y^{\w\tau\rho^{-1}} - y^{\w\tau\rho^{-1}\sigma^{-1}})\w y^{\rho^{-1}})_\kappa \\
& \in & V_\sigma\; . \\
\end{array}
\]
\end{Lemma}

\begin{Lemma}
\label{LemMod9}
There is a commutative diagram of $\b S$-linear maps
\begin{center}
\begin{picture}(450,250)
\put( -50, 200){$\bL/\bLCP$}
\put( 165, 210){\vector(1,0){115}}
\put( 210, 225){$\scm\phi_0$}
\put( 285, 200){${\dis\Ds_{\kappa\in\sClGpp}} \b T_0/V_\kappa$}
\put(  40, 160){\vector(0,-1){110}}
\put(  30, 160){\oval(20,20)[t]}
\put(  10, 105){$\scm\iota$}
\put( 370, 130){\oval(20,20)[t]}
\put( 380, 130){\vector(0,-1){80}}
\put( -20,   0){$\bL/\bLC$}
\put( 150,  10){\vector(1,0){130}}
\put( 210,  25){$\scm\phi$}
\put( 300,   0){${\dis\Ds_{\kappa\in\sClG}} \b T_0/V_\kappa\; $,}
\end{picture}
\end{center}
thus defining $\phi_0$.

\rm
Given $\tau\in G$ and $y\in \b T$, we have $(\tau y)^{q^\mu} = \tau^{q^\mu} y'$ for some $y'\in\b T$. Now $\tau^{q^\mu}$ is a $p'$-element by choice of $\mu$ (\ref{LemOmni} iii).
\end{Lemma}

\begin{Lemma}
\label{LemMod10}
There is an $\b S$-linear map
\[
\begin{array}{rcl}
{\dis\Ds_{\kappa\in\sClGpp}} \b T_0/V_\kappa & \lraa{\psi} & \bL/\bLCP \\
(y_\kappa)_\kappa                            & \lramaps    & \sum_{\kappa\in\sClGpp} \kappa y_\kappa\; . \\
\end{array}
\]

\rm
Given $y\in \Jac(\b T)$ and $\sigma\in\ClGpp$, we have to show that the image of $(\dell_{\sigma,\kappa} y)_\kappa$ vanishes. We have 
$(\sigma y)^{e\cdot o(\sigma)} = (\prod_{i\in [0,o(\sigma) - 1]} y^{\sigma^i})^e = 0$, whence $\sigma y\in\bLCP$.

Now, given $\rho\in C_G(\sigma)$, $y,z\in\b T_0$, we have to show that the images of $(\dell_{\sigma,\kappa} (z^\rho - z))_\kappa$ and
of $(\dell_{\sigma,\kappa}( y^\sigma - y)z)_\kappa$ vanish. Note that 
\[
\sigma y^{\rho^{-1}\sigma} z \; =\; \rho y \sigma z^\rho \rho^{-1}\;\con_{\bLC}\; \sigma y z^\rho\; . 
\]
Putting $\rho = 1$, we see $\sigma y^\sigma z \con_{\bLC} \sigma y z$. Putting $y = 1$, we see $\sigma z \con_{\bLC} \sigma z^\rho$.
\end{Lemma} 

\begin{Lemma}
\label{LemMod10_5}
Let $f\geq 1$. Given $\alpha\in\Gal(\Fu{q^f}/\Fu{q})\iso C_f$ of order $o(\alpha)$ prime to $p$, the map
\[
\begin{array}{rcl}
\Fu{q^f} & \lraa{\sNrm_{\alpha,q^\mu}} & \Fu{q^f} \\
y        & \lramaps                    & \prod_{l\in [0,q^\mu-1]} y^{\alpha^l} \\
\end{array} 
\]
is surjective.

\rm
Let $\alpha = \mb{F}^s$, where $\mb{F}(y) = y^q$. Thus $\Nrm_{\alpha,q^\mu}(y) = y^{\frac{q^{sq^\mu} - 1}{q^s - 1}}$, and so we have to show that $\frac{q^{sq^\mu} - 1}{q^s - 1}$ is 
prime to $q^f - 1$. Since $o(\alpha)$ is assumed to be prime to $p$, we have $f[p] = s[p]$. Writing $f' = f/f[p]$, $s' = s/s[p]$ and
$r = q^{f[p]} = q^{s[p]}$, we have to show that 
\[
\gcd\left(r^{f'} - 1,\frac{r^{s'q^\mu} - 1}{r^{s'} - 1}\right)\; =\; 1\; .
\]
We remark that given $u,v,w\in\Z_{>0}$ with $w\, |\, v$ and $\gcd(w,v/w) = 1$, we have $\gcd(u,v/w) = \gcd(u,v)/\gcd(u,w)$.
Moreover, we remark that given $u,v\in\Z_{> 0}$, we have 
\[
\begin{array}{rcl}
\gcd(r^u - 1, r^v - 1) 
& = & (r-1)\cdot \gcd\left(\frac{r^u - 1}{r-1}, \frac{r^v - 1}{r-1}\right) \\
&\aufgl{Euclid's algorithm} & (r-1)\cdot\frac{r^{\gcd(u,v)} - 1}{r-1} \\
& = &  r^{\gcd(u,v)} - 1\; . \\
\end{array}
\]
Now, since $\gcd(r^{s'} - 1,\frac{r^{s'q^\mu} - 1}{r^{s'} - 1}) = \gcd(r^{s'} - 1,q^\mu) = 1$, the claim follows from
\[
\begin{array}{rcl}
\gcd(r^{f'} - 1,r^{s'q^\mu} - 1)
& = & r^{\gcd(f',s' q^\mu)} - 1 \\
& = & r^{\gcd(f',s')} - 1 \\
& = & \gcd(r^{f'} - 1, r^{s'} - 1)\; . \\
\end{array}
\]
\end{Lemma} 

\begin{Lemma}
\label{LemMod11}
The map $\psi$ is surjective.

\rm
Suppose given $\tau\in G$ and $y\in\b T$. We claim that $\tau y\in\bL/\bLCP$ is contained in the image of $\psi$. By additivity, we may assume that $y$ corresponds to 
$(\dell_{j,i} y)_i$ under $\b T\lraiso \prod_{i\in [1,d]} T/\qfk_i^e$ for some $j\in [1,d]$. Let $e_j$ be the primitive idempotent corresponding to 
$(\dell_{j,i})_i$, so $y = y e_j$.

{\it Case} $e_j^\tau\neq e_j$. Then $e_j^\tau = e_{j'}$ for some $j'\neq j$ since $\tau$ preserves primitive idempotents. Therefore, 
$\tau y = \tau y e_j \con_{\bLC} e_j\tau y = \tau y e_j^\tau = 0$.

{\it Case} $e_j^\tau = e_j$. From $y = y e_j$, we infer $y^{\tau^l} = y^{\tau^l} e_j^{\tau^l} = y^{\tau^l} e_j$ for $l\in\Z$.

We write $\tau = \tau_\rms\tau_\rmu = \tau_\rmu\tau_\rms$ with $o(\tau_\rmu)$ a power of $p$ and $o(\tau_\rms)$ prime to $p$. It suffices to show that there exists a 
$\w y\in\b T$ such that $\tau y \con_{\bLCP} \tau_\rms \w y$. With the ansatz that $\w y$ should correspond to 
$(\dell_{j,i} \w y)_i$ under $\b T\lraiso \prod_{i\in [1,d]} T/\qfk_i^e$, we obtain $(\tau_\rms \w y)^{q^\mu} = \tau^{q^\mu} \prod_{l\in [0,q^\mu - 1]} \w y^{\tau_\srms^l}$,
where $\prod_{l\in [0,q^\mu - 1]}\w y^{\tau_\srms^l}$ corresponds to 
\[
\left(\dell_{j,i} \prod_{l\in [0,q^\mu - 1]} \w y^{\alpha^l}\right)_{\!\! i\in [1,d]}
\]
for some automorphism $\alpha$ of $T/\qfk_j^e$ that fixes $\b S$, of order $o(\alpha)$ prime to $p$. Now $(\tau y)^{q^\mu} = \tau^{q^\mu} y'$, where $y'$ corresponds to 
$(\dell_{j,i} y')_i$. By (\ref{LemMod10_5}), applied to $T/\qfk_j\iso\Fu{q^f}$, the map $T/\qfk_j\lraa{\sNrm_{\alpha,q^\mu}} T/\qfk_j$ is surjective, and so we may 
find $\w y\in\b T$ as in the ansatz such that $(\tau y)^{q^\mu} - (\tau_\rms\w y)^{q^\mu}$ is nilpotent; choose $\mu'\geq 0$ such that the exponent $q^{\mu'}$ annihilates it. Now since 
$\bL/\bLCP\lra \bL/\bLC$, $\xi\lramaps\xi^{q^{\mu + \mu'}}$ is likewise injective (cf.\ \ref{LemOmni} ii, iii), we may conclude from
\[
(\tau y)^{q^{\mu + \mu'}} - (\tau_\rms\w y)^{q^{\mu + \mu'}} \;\auf{\mb{\scr (\ref{LemOmni} i)}}{\con}_{\!\!\!\!\bLC}\; 
\left((\tau y)^{q^\mu} - (\tau_\rms\w y)^{q^\mu}\right)^{q^{\mu'}} \; =\; 0
\]
that $\tau y\con_\bLCP \tau_\rms\w y$.
\end{Lemma}

\begin{Lemma}
\label{LemMod12}
The endomorphism $\psi\phi_0$ is bijective.

\rm
Suppose given an element $(\dell_{\sigma,\kappa} y)$ in $\Ds_{\kappa\in\sClGpp} \b T_0/V_\kappa$.
It is mapped to $\sigma y$ under $\psi$, and then to $(\sigma y)^{q^\mu}$ via $\iota$. Now let $\w\sigma^\rho = \sigma^{q^\mu}$, where $\w\sigma\in\ClGpp$, $\rho\in G$.
Then $(\sigma y)^{q^\mu}$ is mapped via $\phi$ to the tuple that has entry $\left(\prod_{i\in [0,q^\mu - 1]} y^{\sigma^i}\right)^{\rho^{-1}}\in \b T_0/V_{\w\sigma}$ at $\w\sigma$, and
$0$ elsewhere. The corresponding subtuple in $\Ds_{\kappa\in\sClGpp} \b T_0/V_\kappa$ is the image of $(\dell_{\sigma,\kappa} y)$ under $\psi\phi_0$.

Since $\sigma\lramaps\sigma^{q^\mu}$ induces a permutation of $\ClGpp$ and since
\[
\begin{array}{rcl}
\b T_0/V_{\w\sigma} & \lra     & \b T_0/V_{\w\sigma^\rho} \; =\; \b T_0/V_{\sigma^{q^\mu}} \\ 
y                   & \lramaps & y^\rho \\
\end{array}
\]
is bijective, it remains to be shown that 
\[
\begin{array}{rcl}
\b T_0/V_\sigma & \lra       & \b T_0/V_{\sigma^{q^\mu}} \\
y               & \lramaps   & \prod_{i\in [0,q^\mu - 1]} y^{\sigma^i} \\
\end{array}
\]
is bijective. Since $\sigma$ is a $p'$-element of $G$, the subspace $V_{\sigma^{q^\mu}}$ contains the ideal generated by the elements of the form $z^\sigma - z$, where $z\in\b T_0$.
Thus $\prod_{i\in [0,q^\mu - 1]} y^{\sigma^i}\con_{V_{\sigma^{q^\mu}}} y^{q^\mu}$. Moreover, note that $C_G(\sigma^{q^\mu}) = C_G(\sigma)$, and so $V_{\sigma^{q^\mu}} = V_\sigma$.
Since $z\lramaps z^{q^\mu}$ is an automorphism on $\b T_0$, the endomorphism it induces on $\b T_0/V_\sigma$ is an automorphism, too.
\end{Lemma}

\begin{Theorem}
\label{ThBrauer}
We have
\[
\dim_{\b S} Z(\Lambda/\Jac(\Lambda))\; =\; \dim_{\b S} Z(\bL/\Jac(\bL))\; =\; \sum_{\sigma\in\sClGpp}\dim_{\b S}(\b T_0/V_\sigma)\; =:\; \mb{\rm z}(\Lambda)\; .
\]
Concerning notation, cf.\ (\ref{NotMod1}, \ref{NotMod7}).

\rm 
Consider the following diagram of $\b S$-linear maps (cf.\ \ref{LemMod8}, \ref{LemMod9}, \ref{LemMod10}).
\begin{center}
\begin{picture}(900,100)
\put(   0,  40){${\dis\Ds_{\kappa\in\sClGpp}} \b T_0/V_\kappa$}
\put( 280,  50){\vector(1,0){130}}
\put( 330,  65){$\scm\psi$}
\put( 430,  40){$\bL/\bLCP$}
\put( 650,  50){\vector(1,0){130}}
\put( 700,  65){$\scm\phi_0$}
\put( 770,  40){${\dis\Ds_{\kappa\in\sClGpp}} \b T_0/V_\kappa$}
\end{picture}
\end{center}
Since $\psi$ is surjective (\ref{LemMod11}) and $\psi\phi_0$ is bijective (\ref{LemMod12}), $\psi$ is an $\b S$-linear isomorphism. Thus
\[
\sumd{\kappa\in\sClGpp} \dim_{\b S}(\b T_0/V_\kappa) \; =\; \dim_{\b S}(\bL/\bLCP) \;\aufgl{(\ref{LemOmni} iv)}\; \dim_{\b S} Z(\bL/\Jac(\bL))\; .
\]
\end{Theorem}

\begin{Corollary}
\label{CorBrauer0}
The number of isoclasses of simple $\bL$-modules is less than or equal to $\mb{\rm z}(\Lambda)$.
\end{Corollary}

\begin{Corollary}
\label{CorBrauer1}
Suppose $\bL/\Jac(\bL)$ to be split semisimple over $\b S$, i.e.\ suppose all simple $\bL$-modules to be absolutely simple. Then the number of isoclasses of simple $\bL$-modules is 
given by $\mb{\rm z}(\Lambda)$.
\end{Corollary}

We may reformulate (\ref{CorBrauer1}) to an assertion on $T\wr G$.

\begin{Corollary}
\label{CorBrauer2}
Suppose all simple $T\wr G$-modules to be absolutely simple. Then the number of isoclasses of indecomposable projective $T\wr G$-modules is given by $\mb{\rm z}(T\wr G)$.
\end{Corollary}

\begin{Remark}[particular cases]
\label{RemMod13}
We drop the finiteness assumption on $\b S$.
\begin{itemize}
\item[(i)] The $\b T\wr G$-module $\b T_0$ is simple, with endomorphism ring $\End_{\b T\wr G}\b T_0\iso\mb{\rm Fix}_{C_H(\qfk_1)} T/\qfk_1$. If $T/\qfk_1$ is separable over $\b S$, 
then $\b T_0$ is absolutely simple, i.e.\ $\End_{\b T\wr G}\b T_0 = \b S$. Cf.\ (\ref{ExMod3}).
\item[(ii)] If $G$ is a $p$-group, then $\b T_0$ is the only simple $\bL$-module, up to isomorphism. In particular, if $\b S$ is finite, (\ref{ThBrauer}) yields 
$\dim_{\b S} \b T_0/V_1 = \dim_{\b S} \Fix_{\!\! G} \b T_0$.
\item[(iii)] If $\b S$ is finite, and if $G$ acts trivially on $\b T_0$, then $\mb{\rm z}(\Lambda) = |\ClGpp|\cdot \dim_{\b S}\b T_0$ (\ref{ThBrauer}).
\item[(iv)] Suppose $d = 1$. Write $N_0$ for the kernel of the operation of $G$ on $\b T_0$. If $\sigma\in G$ is not contained in $N_0$, then 
$\b T_0/V_\sigma = 0$. Cf.\ (\ref{CorI17}).
\item[(v)] Suppose $G$ to be a $p'$-group. Then $\b T_0\wr G$ is semisimple. In particular, if in addition $T/S$ is unramified, then $\bL = \b T\wr G$ is semisimple. 
\item[(vi)] The isoclasses of simple modules of $\bL = \b T\wr G$ correspond bijectively to the isoclasses of simple modules of $\b T_0\wr G$ via composition of the operation
with $\b T\wr G\lra \b T_0\wr G$. Moreover, we have an isomorphism $(\b T\wr G)/\Jac(\b T\wr G)\lraiso (\b T_0\wr G)/\Jac(\b T_0\wr G)$.
\item[(vii)] If $e = 1$ and if $T/\qfk_1$ is separable over $\b S$, then $V_1$ has codimension $1$ in $\b T_0$.
\item[(viii)] The $\b S$-linear dimension of a $\b T\wr G$-module is divisible by $\dim_{\b S}\b T_0 = h/e$.
\end{itemize}

\rm
Ad (i). Firstly, $\End_{\b T\wr G}\b T_0\iso\mb{Fix}_G \b T_0$. Now the projection of $\mb{Fix}_G \b T_0$ to $T/\qfk_1$ is injective and surjects onto $\mb{Fix}_{C_H(\qfk_1)} T/\qfk_1$.
If the extension $T/\qfk_1$ over $\b S$ is separable, we obtain a surjection of $C_H(\qfk_1)$ onto $\Gal((T/\qfk_1)/\b S)$ {\bf\cite[\rm I.9.4]{N92}}.

Ad (ii). Any $\bL$-module $M$ contains a fixed point different from $0$. In fact, passing to a cyclic module, we may assume it to be finite, of cardinality divisible
by $p$. Since $p$ divides the length of a nontrivial $G$-orbit, the claim follows. Cf.\ {\bf\cite[\rm prop.\ 26]{Se77}}. 

Let $m\in M$ be a nonzero fixed point. Then $m\cdot \b T$ is a submodule of $M$. Hence, if $M$ is simple, $M$ is isomorphic to a quotient of $\b T$. But the only simple quotient of
$\b T$ is $\b T_0$.

Ad (iv). By assumption, there exists a $y\in\b T_0$ such that $y^\sigma - y \neq 0$. The subspace $V_\sigma$ of $\b T_0$ contains the ideal generated by $y^\sigma - y$. Since
now $\b T_0$ is a field, this implies $V_\sigma = \b T_0$.

Ad (v). Given an epimorphism of $\bL$-modules, we may choose a $\b T_0$-linear coretraction, for $\b T_0$ is semisimple. The sum over the $G$-conjugates divided by $|G|$ yields
a $\bL$-linear coretraction.

Ad (vi). The kernel of $\b T\wr G\lra \b T_0\wr G$ is a nilpotent ideal, thus contained in $\Jac(\b T\wr G)$.

Ad (vii). By tameness of $T/S$, (\ref{RemI3_2}) and {\bf\cite[\rm 32.1]{CR81}} show that $\b T_0|_{\b S H}\iso\b S H$. Under such an isomorphism, $V_1$ corresponds to the 
augmentation ideal.

Ad (viii). We are reduced to consider a simple $\b T\wr G$-module $X$. Let $e_j$ correspond to  $(\dell_{j,i})_i$ under $\b T_0\lraiso\prod_{i\in [1,d]} T/\qfk_i$. We may decompose
$X$ as a $\b T_0$-module into $X = \Ds_{j\in [1,d]} X e_j$. Since $G$ acts transitively on $\{ e_j\; |\; j\in [1,d]\}$, we have for each $j\in [1,d]$ an element $\sigma\in G$ that 
induces $X e_1\lraiso X e_j$ by multiplication. Since $X e_1$ is a module over $T/\qfk_1$, $d\cdot\dim_{\b S} T/\qfk_1 = \dim_{\b S_0}\b T_0$ divides $\dim_{\b S} X$.
\end{Remark}

\subsection{Examples}

\begin{quote}
\begin{footnotesize}

When we count simple modules, we shall tacitly count them up to isomorphism. We use {\sc Maple,} {\sc Magma} and the {\sc Meat Axe} version of {\sc M.\ Ringe} {\bf\cite{Ri93}}.

\begin{Example}
\label{ExMod1}\rm
We continue (\ref{ExI2_6}).
\begin{itemize}
\item[$p = 2$.] Let $S = \Z_{(2)}$ and $T = \Z_{(2)}[i]$. Since $\b T_0 \iso \Fu{2}$, the simple $\b T\wr G$-modules are given by the two simple $\Fu{2}\Sl_3$-modules 
(\ref{RemMod13} iii).
\item[$p = 3$.] Let $S = \Z_{(3)}$ and $T = \Z_{(3)}[i]$. Now $\b T = \b T_0 = \Fu{3}[i] \iso \Fu{9}$. We obtain $V_1 = \Fu{3}\spi{i}$ and $V_{(1,2)} = \Fu{3}[i]$. 
Therefore $\mb{z}(\Z_{(3)}[i]\wr \Sl_3) = 1$, and thus $\b T_0$ is the unique (absolutely) simple $\b T\wr G$-module (\ref{RemMod13} i). 
\end{itemize}
\end{Example}

\begin{Example}
\label{ExMod2}\rm
Let $K = \Q$. Let $\gamma$ be a root of the irreducible polynomial 
\[
\mu(X) = X^6 - 3 X^5 + 7 X^4 - 9 X^3 + 7 X^2 - 3 X + 1\; \in\;\Q[X]\; ,
\]
and let $L = \Q(\gamma)$. Then $\Z[\gamma]$ is the ring of algebraic integers in the Galois extension $L/K$. We have $\Gal(L/K) \iso \Sl_3$. Thus there is an
operation of $\Sl_4$ on $L$, where
\[
\begin{array}{rcl}
\gamma^{(1,2)}     & = & 1 - \gamma \\
\gamma^{(1,2,3,4)} & = & 1/\gamma\; . \\
\end{array}
\]

\begin{itemize}
\item[$p = 2$.] Let $S = \Z_{(2)}$ and $T = \Z_{(2)}[\gamma]$. We have $\mu(X)\con_2 (X^3 + X + 1)(X^3 + X^2 + 1)$, hence $\b T = \b T_0 \iso \Fu{8}\ti \Fu{8}$. Let
$\mb{F}x = x^2$ for $x\in\Fu{8}$. The operation of $\Sl_4$ on $\b T_0$ is given by
\[
\begin{array}{rclcrcl}
\Fu{8} & \ti & \Fu{8} & \lra                 & \Fu{8}  & \ti & \Fu{8} \\
x      & \ti & y      & \lramapsa{(1,2)}     & y       & \ti & x \\
x      & \ti & y      & \lramapsa{(1,2,3,4)} & \mb{F}y & \ti & \mb{F}^2 x\; . \\
\end{array}
\]
Let $\alpha\in\Fu{8}$ with $\alpha^3 + \alpha + 1 = 0$. Then $V_1 = \Fu{2}\spi{1\ti 1, \alpha\ti 0, \alpha^2\ti 0, 0\ti\alpha, 0\ti\alpha^2}$ (cf.\ \ref{RemMod13} vii), and 
$V_{(1,2,3)} = \b T_0$. Altogether, $\mb{z}(T\wr G) = 1$. Therefore, $\b T_0 = \b T$ is the only (absolutely) simple $\b T\wr G$-module (\ref{ThBrauer}, \ref{RemMod13} i).

\item[$p = 3$.] Let $S = \Z_{(3)}$ and $T = \Z_{(3)}[\gamma]$. We have $\mu(X)\con_3 (X^2 + 1)(X^2 + X - 1)(X^2 - X - 1)$, hence 
$\b T = \b T_0 \iso \Fu{9}\ti \Fu{9}\ti \Fu{9}$. Let $\mb{F}x = x^3$ for $x\in\Fu{9}$. The operation of $\Sl_4$ on $\b T_0$ is given by
\[
\begin{array}{rccclcrcccl}
\Fu{9} & \ti & \Fu{9} & \ti & \Fu{9} & \lra                 & \Fu{9}  & \ti & \Fu{9} & \ti & \Fu{9} \\
x      & \ti & y      & \ti & z      & \lramapsa{(1,2)}     & y       & \ti & x      & \ti & \mb{F}z \\
x      & \ti & y      & \ti & z      & \lramapsa{(1,2,3,4)} & \mb{F}x & \ti & z      & \ti & y \; . \\
\end{array}
\]
Choose $\iota\in\Fu{9}$ with $\iota^2 + 1 = 0$. Then $V_1 = \Fu{3}\spi{1\ti -1\ti 0, 0\ti 1\ti -1, \iota\ti 0\ti 0, 0\ti\iota\ti 0, 0\ti 0\ti\iota}$ (cf.\ \ref{RemMod13} vii).
Moreover, $V_{(1,2)(3,4)} = \Fu{3}\spi{\iota\ti 0\ti 0, 0\ti 1\ti -1, 0\ti\iota\ti -\iota}$. Finally, $V_{(1,2)} = V_{(1,2,3,4)} = \b T_0$. Altogether, $\mb{z}(T\wr G) = 4$.
It turns out that there are three simple $\b T\wr G$-modules, two of dimension $6$ and endomorphism ring $\Fu{3}$, and one of dimension $12$ and endomorphism ring $\Fu{9}$.
In particular, the center of $\b T\wr G/\Jac(\b T\wr G)$ is isomorphic to $\Fu{3}\ti \Fu{3}\ti \Fu{9}$. It is of dimension $4$ over $\Fu{3}$, as predicted by (\ref{ThBrauer}).

\item[$p = 31$.] Let $S = \Z_{(31)}$ and $T = \Z_{(31)}[\gamma]$. We have $\mu(X)\con_{31} (X-2)^2 (X + 15)^2 (X+1)^2$, hence 
$\b T \iso \Fu{31}[\eps]\ti\Fu{31}[\eps]\ti\Fu{31}[\eps]$, where $\eps^2 = 0$. In particular, $\b T_0 \iso \Fu{31}\ti\Fu{31}\ti\Fu{31}$. The operation of $\Sl_4$ on $\b T_0$ is given by
\[
\begin{array}{rccclcrcccl}
\Fu{31} & \ti & \Fu{31} & \ti & \Fu{31} & \lra                 & \Fu{31} & \ti & \Fu{31} & \ti & \Fu{31} \\
x       & \ti & y       & \ti & z       & \lramapsa{(1,2)}     & z       & \ti & y       & \ti & x \\
x       & \ti & y       & \ti & z       & \lramapsa{(1,2,3,4)} & y       & \ti & x       & \ti & z \; . \\
\end{array}
\]
We obtain $V_1 = \Fu{31}\spi{1\ti -1\ti 0, 1\ti 0\ti -1}$, $V_{(1,2)} = \Fu{31}\ti 0\ti \Fu{31}$, $V_{(1,2)(3,4)} = \Fu{31}\spi{1\ti 0\ti -1}$, $V_{(1,2,3)} = \b T_0$ and
$V_{(1,2,3,4)} = \Fu{31}\ti\Fu{31}\ti 0$. Altogether, $\mb{z}(T\wr G) = 5$. In fact, there are $5$ absolutely simple $\b T\wr G$-modules, four of dimension $3$ and one of dimension $6$.
\end{itemize}
\end{Example}

\begin{Example}
\label{ExMod3}\rm
Let $G = H = C_2 = \spi{c}$, let $L = \Q(X)$, let $X^c = -X$, so that $K = \Q(X^2)$. Let $\pfk := 2\Z[X^2]\tm\Z[X^2]$, let $S = \Z[X^2]_\pfk$, hence $T = \Z[X]_\pfk$. 
We obtain $\b S = \Fu{2}(X^2)$ and $\b T = \b T_0 = \Fu{2}(X)$. In particular, $c$ acts trivially on $\b T_0$, so $\End_{\b T\wr G}\b T_0 = \Fu{2}(X)$. Therefore, $\b T_0$ is
simple, but not absolutely simple over the $\b S$-algebra $\b T\wr G$.
\end{Example}
\end{footnotesize}
\end{quote}

\parskip0.0ex
\begin{footnotesize}

\parskip1.2ex

\vspace*{1cm}

\begin{flushright}
Matthias K\"unzer\\
Universit\"at Ulm\\
Abteilung Reine Mathematik\\
D-89069 Ulm\\
kuenzer@mathematik.uni-ulm.de\\
\end{flushright}
\end{footnotesize}
\end{document}